\documentclass[a4paper,12pt]{article}
\usepackage{graphicx}
\usepackage{amsmath}
\usepackage{amsthm}
\usepackage{amssymb}
\usepackage{amsfonts}
\usepackage{latexsym}
\usepackage{sidecap}
\usepackage{CJK,graphicx}
\usepackage{framed}
\usepackage{enumerate}
\usepackage{float}
\newcommand{\new}{\newcommand*}
\new{\pt}{\partial}
\new{\nn}{\nonumber}
\new{\R}{{\sf R}}
\new{\RS}{{\sf S}}
\new{\bk}{\mathbf{k}}
\new{\br}{\mathbf{r}}
\new{\bu}{\boldsymbol{u}}
\new{\bn}{\boldsymbol{n}}
\new{\newt}{\newtheorem}
\new{\sect}[1]{\section{#1}\setcounter{equation}{0}}

 \title{High efficiently numerical simulation of the TDGL equation with reticular free energy in hydrogel}
 \author{ Jun Han  \and Hui Zhang\and Zhengru Zhang \thanks{School of Mathematical Sciences, Beijing Normal University, Beijing, 100875, P.R. China .}}

\begin{document}

\maketitle
\begin{abstract}
 In this paper, we focus on  the numerical simulation of phase separation about macromolecule microsphere composite (MMC) hydrogel. The model equation is based on Time-Dependent Ginzburg-Landau (TDGL) equation with reticular free energy. We have put forward two $L^2$ stable schemes to simulate simplified TDGL equation. In numerical experiments, we observe that simulating the whole process of phase separation requires a considerably long time. We also notice that the total free energy changes significantly in initial time and varies slightly in the following time. Based on these properties, we introduce an adaptive strategy based on one of stable scheme mentioned. It is found that the introduction of the time adaptivity cannot only resolve the dynamical changes of the solution accurately but also can significantly save CPU time for the long time simulation.\\
\end{abstract}

{\bf Keywords:} TDGL equation, $L^2$-norm stale, semi-implicit finite difference schemes, adaptive time-stepping strategy \\

\section{Introduction}
Hydrogels, which are polymer with 3-D crosslinked hydrophilic network structure, have increasingly extensive applications in industrial and biomedical fields [8,9,22]. Therefore, hydrogels have recently received tremendous attention in scientific communities due to their enhanced properties [1,15]. Traditional hydrogels usually have poor mechanism as they consist of water in a fragile polymer network and the cross-linking points in the system are disordered leading to different chain length, which causes shorter chains bearing more stress crack and the whole hydrogel breaks down. However, in recent year, there has been a lot of novel methods to improve the structure of hydrogels and thus significantly enhancing their mechanical properties [8,11]. These methods primarily introduce special corss-linking center in gel structure such as sliding crosslinking and multiple crosslinking center, which can effectively disperse stress and thus the polymer is not easy to fracture [9].   According to the structure of these novel hydrogels, they can be classified as topological (TP) gels, nanocomposite (NC) hydrogels, double-network (DN) gels, MMC hydrogels and so on [8,9]. TP gels have eight crosslinkers that can slide along polymer chains [3] and can swells to about 500 times its original weight and can be stretched to nearly 20 times its original length [14]. NC hydrogels, which utilize inorganic material to form crosslinking junctions and adsorb strongly on  surfaces of platelets with enough chains attached to provide bridges, constitute a network structure jielding unusual mechanical properties including very good toughness [7]. DN gels generally consist of two independently cross-linked networks, one consisting of a rigid polyelectrolyte and the other a flexible uncharged polymer, and exhibit extremely high mechanical strength [4,6]. MMC hydrogels, which were synthesized by Ting Huang in 2007 [8], have a unigue well-defined microstructure and very high mechanical strength. MMC hydrogels are environmentally sensitive and mainly possible to use in drug delivery and other biomedical applications [20]. Much work in chemical structure and dynamics simulation have been done about TP gels, NC hydrogels, and DN gels and much progress has been made [5,21,10,11]. But there are relatively less results about MMC hydrogels.

Meanwhile, many computational methods have been developed and applied in modeling and simulation of phase separation of hydrogels. They range from molecular scale (e.g., molecular dynamics, Monte Carlo), microscale (e.g., Brownian dynamics, dissipative particle dynamics, lattice Boltzmann, TDGL method, dynamic density functional theory method) to mesoscale and macroscale (e.g., micromechanics, equivalent-continuum and self-similar approaches, finite element method) [21]. These methods are foundation to simulate phase separation of hydrogels. However, observing the whole process of phase separation is computationally expensive and costs a large amount of time when these methods are used. Therefore, seeking for algorithms to decrease computation cost is an urgent task. The adaptive time-stepping method is introduced here. Adaptive time stepping has been well studied for solving imitial value problems in ODEs. The review in [17,18] summarizes several step control methods for local time adaptivity based on linear feedback theory. A locally varying time step method is developed for solving hyperbolic conservative PDEs where large time step is adopted on region with smooth solution while small time step is taken in domain with nearly singular solution within same time level [19]. In [16], the adaptive time-stepping technique has been developed based on energy variation which is an important physical quantity in the molecular beam epitaxy growth model. Following [16], we improve the choice of adaptive time step according to properties of TDGL equation.

In this paper, considering that MMC hydrogels have microstructure with network crosslinking points, we exploit TDGL method to simulate phase evolution of MMC hydrogels based on reticular free energy deduced in [22]. TDGL method is an accurate and robust method to simulate structural evolution of phase separation in polymer blends and copolymers and many nice results have been gotten in numerical experiments [1,13,21]. There are two main objectives in this paper. First, we put forward two $L^2$ stable second-order implicit schemes to simulate the equation, which satisfy the property of total free energy decreasing. One scheme needs to solve a high-dimensional matrix in each step, which can be solved by conjugate gradient method. The other is a nonlinear scheme where we need Newton iteration in each step, which will cost more computation. Second, as the equation is unconditionally energy stale and the schemes are $L^2$ norm stable, we introduce an adaptive strategy which is vert effective for such evolutionary equations where total free energy decays significantly in initial time and slightly later, called "rough-smooth" pattern. During the time with rapid free energy decay, small adaptive time steps are chosen adaptively, while large time steps are used when free energy is varied smoothly.

This paper is presented as follows. In section 2, TDGL equation is introduced in detail with reticular free energy. Then, we employ variational principle to derive and simplify the equation in a further form, which will be used in following theoretical analysis and numerical simulation. In section 3, analysis of the equation is given and two semi-implicit schemes are proposed.
In section 4, some numerical experiments are conducted to test one of the numerical schemes. In section 5, based on observations, we put forward a suitable adaptive time-stepping strategy and some concluding remarks are drawn in final section.\\

\section{TDGL equation}
TDGL equation is a microscale method for simulating the structural evolution of phase-separation in polymer blends and block copolymers. For an incompressible polymer blend, the model equation can be written as follows [1,13,21,22].
\begin{equation}\label{eq1}
\frac{\partial\phi(\textbf{x},t)}{\partial t}=\nabla\cdot(M(\phi)\nabla\frac{\delta U[\phi(\textbf{x},t)]}{\delta \phi(\textbf{x},t)})+\zeta(\textbf{x},t),
\end{equation}
where $\phi(\textbf{x},t)$ is order parameter, by definition $\phi \in(0,1)$, the mobility $M$ may depend on $\phi$. $\textbf{x} \in{\Omega}$, $\Omega =(0,2\pi)\times(0,2\pi)$. The total Flory-Huggins-de Gennes free energy $U$ can be defined as:
\begin{equation}\label{eq2}
U[\phi(\textbf{x},t)]=k_BT\int_\Omega[\frac{F(\phi(\textbf{x},t))}{k_BT}+\kappa(\phi(\textbf{x},t))|\nabla\phi(\textbf{x},t)|^2]d\textbf{x},
\end{equation}
where $\kappa(\phi)=\frac{\sigma^2}{36\phi(1-\phi)}$, $k_B$ is Boltzmann constant, $\sigma$ is the Huhn length of the polymer,$T$ is temperature constant [2]. It is pointed out that for a binary polymer blend, $F(\phi(\textbf{x},t))$ usually denotes the Flory-Huggins free energy. However, as mentioned before, MMC hydrogels have network structure. Hence, in [8], it is proven that the reticular free energy should be defined as follows.
\begin{equation}\label{eq3}
F=k_BT[\frac{\phi}{\tau}ln(\frac{\alpha\phi}{\tau})+\frac{\phi}{N}ln(\frac{\beta\phi}{\tau})+(1-\rho\phi)ln(1-\rho\phi)+
\chi\rho\phi(1-\rho\phi)],
\end{equation}
where parameters $\tau$, $\alpha$, $\beta$, $N$, and $\rho$ are constant variables [22]. $\chi$ is the enthalpic interaction parameter between two polymer components. As derived in [22], it can be guaranteed that $(1-\rho\phi)>0$.

The thermal noise in (\ref{eq1}) is a random term with zero mean which should satisfy following conditions:
$$
<\zeta(\textbf{x},t),\zeta(\textbf{x}',t)>=\epsilon M_0\nabla^2\delta(\textbf{x}-\textbf{x}')\delta(t-t'),
$$
where $\epsilon$ is the small magnitude of the fluctuation [12].\\

\subsection{Analysis of energy stability}
When the property of energy stability is considered, the noise term in (\ref{eq1}) should not be considered. Therefore, multiplying (\ref{eq1}) (deleting noise term) by $\frac{\delta U}{\delta \phi}$ and integrating over $\Omega$, we get
$$
\frac{dU}{dt} =\frac{1}{|\Omega|}\int_\Omega\frac{\delta U(\phi)}{\delta \phi}\nabla\cdot(M(\phi)\nabla\frac{\delta U(\phi)}{\delta \phi})d\textbf{x}.
$$

Here, $|\Omega|$ is the volume of $\Omega$. The above equality holds because of the fact that $U$ is independent of variable $\textbf{x}$ according to definition of $U$ in (\ref{eq2}). Then we can easily derive
$$
\frac{dU}{dt}=\frac{1}{|\Omega|}(\int_\Omega\nabla\cdot(\frac{\delta U(\phi)}{\delta \phi}M(\phi)\nabla\frac{\delta U(\phi)}{\delta \phi})d\textbf{x}-
\int_\Omega M(\phi)\nabla\frac{\delta U(\phi)}{\delta \phi}\cdot \nabla\frac{\delta U(\phi)}{\delta \phi}d\textbf{x}).
$$
By divergence theorem and considering periodic boundary condition for $\phi$, the first term on the right side of above equality disappears. Hence, we get
\begin{equation}\label{eq4}
\frac{dU}{dt}=-\frac{1}{|\Omega|}\int_\Omega M(\phi)|\nabla\frac{\delta U(\phi)}{\delta \phi}|^2d\textbf{x}\leq 0.
\end{equation}
Therefore, we know that the total free energy of whole system is decreasing over time, i.e, the system is energy stable.\\

\subsection{Deduction of the model equation}
For any $\psi\in C^1(\Omega)$ and considering periodic boundary condition for $\phi$ and $\psi$, we define
$$
H(\epsilon)=U[\phi+\epsilon\psi]=k_BT\int_\Omega[\frac{F(\phi(\textbf{x})+\epsilon\psi(\textbf{x}))}{k_BT}+\kappa(\phi(\textbf{x})+
\epsilon\psi(\textbf{x}))|\nabla(\phi(\textbf{x})+\epsilon\psi(\textbf{x}))|^2]d\textbf{x}.
$$
Obviously, $H(\epsilon)\in C^1(\mathbb{R})$. Hence, $H'(\epsilon)$ can be easily derived and when $\epsilon=0$, we get
\begin{equation}\label{eq5}
H'(0)=k_BT\int_\Omega[\frac{F'(\phi)\psi}{k_BT}+\kappa'(\phi)\psi|\nabla\phi|^2+
2\kappa(\phi)\nabla\phi\cdot\nabla\psi]d\textbf{x}.
\end{equation}
According to the Green's formula and periodic boundary condition for $\phi$ and $\psi$, we derive
$$
\int_\Omega\kappa(\phi)\nabla\phi\cdot\nabla\psi d\textbf{x}=\int_{\Omega}\nabla\cdot(\kappa(\phi)\psi\nabla\phi)d\textbf{x}
-\int_\Omega\nabla\cdot(\kappa(\phi)\nabla\phi)\psi d\textbf{x}
$$
$$
=-\int_\Omega[\kappa'(\phi)\nabla\phi\cdot\nabla\phi+\kappa(\phi)\Delta\phi]\psi d\textbf{x}.
$$
By the variational principle in (\ref{eq5}), we get
\begin{equation}\label{eq6}
\frac{\delta U[\phi]}{\delta \phi}=H'(0)=F'(\phi)
-k_BT(\kappa'(\phi)|\nabla\phi|^2+2\kappa(\phi)\Delta\phi).
\end{equation}
Combining (\ref{eq1}) with (\ref{eq6}), following [1,13,21,22] where the mobility $M$ is equal to constant $M_0$, we get the model equation
\begin{equation}\label{eq7}
\frac{\partial \phi}{\partial t}=M_0\nabla\cdot(\nabla F'(\phi))-M_0k_BT\Delta(\kappa'(\phi)|\nabla\phi|^2)-2M_0k_BT\Delta(\kappa(\phi)\Delta\phi)+\zeta(\textbf{x},t).
\end{equation}
Straightforward computations in (\ref{eq3}) jield $\nabla F'(\phi)=k_BT(\frac{1}{\tau\phi}+\frac{1}{N\phi}+\frac{\rho^2}{1-\rho\phi}-2\chi\rho^2)\nabla\phi$. Here, we denotes $G(\phi)=k_BT(\frac{1}{\tau\phi}+\frac{1}{N\phi}+\frac{\rho^2}{1-\rho\phi}-2\chi\rho^2)$.

In practical problems in materials, many chemists and physicists have taken $\kappa(\phi)$ to be a constant parameter. Therefore, in practical simulation, it is reasonable to assume that $\kappa(\phi)$ is constant, denoted by $K$. Then, the model equation would be simplified as:
\begin{equation}\label{eq8}
\frac{\partial \phi}{\partial t}=M_0G'(\phi)|\nabla \phi|^2+M_0G(\phi)\Delta\phi-2M_0k_BTK\Delta^2\phi+\zeta(\textbf{x},t).
\end{equation}

\section{$L^2$ stable semi-implicit finite difference schemes}
As illustrated in [23], A system located at a metastable state initially will fluctuate around this metastable state under thermal noise perturbation. When the noise happens to be large enough to carry the system out of the basin of attraction of this metastable state, the phase separation will occur. The role of noise therm is to carry the system out of the basin of the metastable state. Therefore, in theoretical analysis, we do not need to consider noise term and the model equation can be rewritten:
\begin{equation}\label{eq9}
\frac{\partial \phi}{\partial t}=M_0G'(\phi)|\nabla \phi|^2+M_0G(\phi)\Delta\phi-2M_0k_BTK\Delta^2\phi.
\end{equation}

Before we construct finite difference schemes to solve equation (\ref{eq8}), we first prove that the analytic solution of (\ref{eq8}) is $L^2$ stable if it existed, which can be rewritten in following theorem.

THEOREM 3.1. If $\phi(\textbf{x},t)$ is a solution of Eq. (\ref{eq9}), then $\phi(\textbf{x},t)$ is $L^2$ stable with respect to the initial condition:
\begin{equation}\label{eq10}
\|\phi(\cdot,t)\|^2\leq\|\phi(\cdot,0)\|^2,
\end{equation}
where $\|\cdot\|$ is the standard $L^2$-norm in $\Omega$.

Proof: Multiplying (\ref{eq9}) by $\phi$ and integrating over $\Omega$, we have
$$
\frac{1}{2}\frac{d}{dt}\int_\Omega|\phi|^2d\textbf{x}=M_0\int_\Omega G'(\phi)\phi|\nabla\phi|^2d\textbf{x}+M_0\int_\Omega G(\phi)\phi\Delta\phi d\textbf{x}-2M_0k_BTK\int_\Omega\phi\Delta^2\phi d\textbf{x}.
$$
By the Green's formula for the second term and third term on the right side of the above equality and considering the periodic boundary condition for $\phi$, we get
$$
\frac{1}{2}\frac{d}{dt}\int_\Omega|\phi|^2d\textbf{x}=M_0\int_\Omega G'(\phi)\phi|\nabla\phi|^2 d\textbf{x}-M_0\int_\Omega \nabla (G(\phi)\phi)\cdot\nabla\phi d\textbf{x}-2M_0k_BTK\int_\Omega |\Delta\phi|^2d\textbf{x}.
$$
After simple calculations, we derive
\begin{equation}\label{eq11}
\frac{1}{2}\frac{d}{dt}\int_\Omega|\phi|^2d\textbf{x}=-M_0\int_\Omega G(\phi)|\nabla\phi|^2d\textbf{x}
-2M_0k_BTK\int_\Omega |\Delta\phi|^2d\textbf{x}.
\end{equation}

According to the value of parameter $\rho,\chi,N$, and $\tau$, it is obvious that $G(\phi)\geq0$, Therefore, the following inequality holds:
$$
\frac{d}{dt}\int_\Omega|\phi|^2d\textbf{x}\leq0.
$$
Hence, inequality (\ref{eq10}) is proven. \\

\subsection{A second-order semi-implicit linear scheme}
In the following, we construct finite difference schemes to solve the Eq. (\ref{eq9}). First, we give a semi-implicit linear scheme with second-order as follows.
$$
\frac{\phi_{i,j}^{n+1}-\phi_{i,j}^{n}}{\Delta t}=
M_0G'(\phi_{i,j}^n)\frac{|\nabla_h\phi_{i,j}^{n}|^2+\nabla_h\phi_{i,j}^{n}\cdot\nabla_h\phi_{i,j}^{n+1}}{2}+M_0G(\phi_{i,j}^n)\frac{\Delta_h\phi_{i,j}^{n}+
\Delta_h\phi_{i,j}^{n+1}}{2}
$$
\begin{equation}\label{eq12}
-2M_0k_BTK\frac{\Delta_h^2\phi_{i,j}^{n+1}+\Delta_h^2 \phi_{i,j}^{n}}{2}.
\end{equation}
$\nabla_h$,$\Delta_h$ are central-scheme operators of finite difference method, and h is the mesh size of the spatial discretizition. Obviously, in each step, (\ref{eq11}) can be solved by computing a linear system of equations. Thus, we call it as a linear scheme.

THEOREM 3.2. let ${\phi_{i,j}^n}$ be the numerical solution of (\ref{eq9}) at $t_n$. Then we have
\begin{equation}\label{eq13}
\parallel\phi^{n}\parallel_h^2\leq\parallel\phi^{0}\parallel_h^2,
\end{equation}
where $\parallel\cdot\parallel_h$ is the discrete $L^2$-norm.

Proof. Multiplying (\ref{eq12}) with the rewritten form by $\frac{\phi_{i,j}^{n+1}+\phi_{i,j}^{n}}{2}$, we get
$$
\frac{(\phi_{ij}^{n+1})^2-(\phi_{ij}^{n})^2}{2\Delta t}=
M_0G'(\phi_{i,j}^n)\nabla_h\phi_{i,j}^{n}\frac{\nabla_h\phi_{i,j}^{n}+\nabla_h\phi_{i,j}^{n+1}}{2}\frac{\phi_{i,j}^{n+1}+\phi_{i,j}^{n}}{2}
$$
$$
+M_0G(\phi_{i,j}^n)\frac{\phi_{i,j}^{n+1}+\phi_{i,j}^{n}}{2}\frac{\Delta_h\phi_{i,j}^{n}+
\Delta_h\phi_{i,j}^{n+1}}{2}-2M_0k_BTK\frac{\Delta_h^2\phi_{i,j}^{n+1}+
\Delta_h^2\phi_{i,j}^{n}}{2}\cdot\frac{\phi_{i,j}^{n+1}+\phi_{i,j}^{n}}{2}.
$$
Summing for $i=1,2,3,...,N_x$ and $j=1,2,...,N_y$. By the Green's formula and considering periodic boundary condition for $\phi$, the above equality can be rewritten:
$$
\frac{1}{2\Delta t}(\parallel\phi^{n+1}\parallel_h^2-\parallel\phi^{n}\parallel_h^2)=
(M_0G'(\phi^n)\nabla_h\phi^{n}\frac{\nabla_h\phi^{n}+\nabla_h\phi^{n+1}}{2},\frac{\phi^{n+1}+\phi^{n}}{2})_h+
$$
$$
(M_0G(\phi^n)\frac{\phi^{n+1}+\phi^{n}}{2},\frac{\Delta_h\phi^{n}+
\Delta_h\phi^{n+1}}{2})_h
-2M_0k_BTK\parallel\frac{\Delta_h(\phi^{n+1}+\phi^{n})}{2}\parallel_h^2.
$$
where $(\cdot,\cdot)_h$ denotes the discrete $L^2$ inner product. By Green's formula for second term on right side above and considering periodic boundary condition for $\phi$, the above equality becomes
$$
\frac{1}{2\Delta t}(\parallel\phi^{n+1}\parallel_h^2-\parallel\phi^{n}\parallel_h^2)=
(M_0G'(\phi^n)\nabla_h\phi^{n}\frac{\nabla_h\phi^{n}+\nabla_h\phi^{n+1}}{2},\frac{\phi^{n+1}+\phi^{n}}{2})_h
$$
$$
-(M_0\nabla_h(G(\phi^n)\frac{\phi^{n+1}+\phi^{n}}{2}),\frac{\nabla_h\phi^{n}+
\nabla_h\phi^{n+1}}{2})_h
-2M_0k_BTK\parallel\frac{\Delta_h(\phi^{n+1}+\phi^{n})}{2}\parallel_h^2.
$$
The second term on right side above can be decomposed to two terms, one of them can be varnished by the first term on right side. Thus, we can get
$$
\frac{1}{2\Delta t}(\parallel\phi^{n+1}\parallel_h^2-\parallel\phi^{n}\parallel_h^2)=
-(M_0G(\phi^n),\frac{|\nabla_h\phi^{n+1}+\nabla_h\phi^{n}|^2}{2})_h
$$
$$
-2M_0k_BTK\parallel\frac{\Delta_h(\phi^{n+1}+\phi^{n})}{2}\parallel_h^2.
$$
Consequently, the right side above is less than 0, we get
$$
\frac{1}{2\Delta t}(\parallel\phi^{n+1}\parallel_h^2-\parallel\phi^{n}\parallel_h^2)\leq0.
$$
Hence, inequality (\ref{eq13}) holds.\\

\subsection{A semi-implicit nonlinear scheme with second-order}
Similar to the linear scheme, we put forward a nonlinear scheme where the equation should be solved by the Newton's iteration.
$$
\frac{\phi_{i,j}^{n+1}-\phi_{i,j}^{n}}{\Delta t}=
M_0G'(\frac{\phi_{i,j}^{n+1}+\phi_{i,j}^{n}}{2})|\nabla_h\frac{\phi_{i,j}^{n}+\phi_{i,j}^{n+1}}{2}|^2+M_0G(\frac{\phi_{i,j}^{n+1}+\phi_{i,j}^{n}}{2})
\frac{\Delta_h\phi_{i,j}^{n}+\Delta_h\phi_{i,j}^{n+1}}{2}
$$
\begin{equation}\label{eq14}
-2M_0k_BTK\frac{\Delta_h^2\phi_{i,j}^{n+1}+\Delta_h^2 \phi_{i,j}^{n}}{2}.
\end{equation}

THEOREM 3.3. let ${\phi_{i,j}^n}$ be the numerical solution of (\ref{eq14}) at $t_n$. Then we have
\begin{equation}\label{eq15}
\parallel\phi^{n}\parallel_h^2\leq\parallel\phi^{0}\parallel_h^2.
\end{equation}

Proof. Similar to the proof of Th.3.2, multiply (\ref{eq14}) by $\frac{\phi_{i,j}^{n+1}+\phi_{i,j}^{n}}{2}$ and summing for $i=1,2,3,...,N_x$ and $j=1,2,...,N_y$. By the Green's formula and periodic boundary condition for $\phi$ is considered, we get
$$
\frac{1}{2\Delta t}(\parallel\phi^{n+1}\parallel_h^2-\parallel\phi^{n}\parallel_h^2)=
(M_0G'(\frac{\phi^n+\phi^{n+1}}{2})|\frac{\nabla_h\phi^{n}+\nabla_h\phi^{n+1}}{2}|^2,\frac{\phi^{n+1}+\phi^{n}}{2})_h+
$$
$$
(M_0G(\frac{\phi^n+\phi^{n+1}}{2})\frac{\phi^{n+1}+\phi^{n}}{2},\frac{\Delta_h\phi^{n}+
\Delta_h\phi^{n+1}}{2})_h
-2M_0k_BTK\parallel\frac{\Delta_h(\phi^{n+1}+\phi^{n})}{2}\parallel_h^2.
$$
By the Green's formula for the second term on right side above and considering the periodic boundary condition for $\phi$, the above equality becomes
$$
\frac{1}{2\Delta t}(\parallel\phi^{n+1}\parallel_h^2-\parallel\phi^{n}\parallel_h^2)=
(M_0G'(\frac{\phi^n+\phi^{n+1}}{2})|\frac{\nabla_h\phi^{n}+\nabla_h\phi^{n+1}}{2}|^2,\frac{\phi^{n+1}+\phi^{n}}{2})_h
$$
$$
-(M_0\nabla_h(G(\frac{\phi^n+\phi^{n+1}}{2})\frac{\phi^{n+1}+\phi^{n}}{2}),\frac{\nabla_h\phi^{n}+
\nabla_h\phi^{n+1}}{2})_h
-2M_0k_BTK\parallel\frac{\Delta_h(\phi^{n+1}+\phi^{n})}{2}\parallel_h^2.
$$
Similarly, the second term on the right side above can be decomposed to two terms, one of which is varnished by first term right side. Therefore, we derive
$$
\frac{1}{2\Delta t}(\parallel\phi^{n+1}\parallel_h^2-\parallel\phi^{n}\parallel_h^2)=
-(M_0G(\frac{\phi^{n+1}+\phi^{n}}{2}),
\frac{|\nabla_h\phi^{n+1}+\nabla_h\phi^{n}|^2}{2})_h
$$
$$
-2M_0k_BTK\parallel\frac{\Delta_h(\phi^{n+1}+\phi^{n})}{2}\parallel_h^2.
$$
Therefore, for $\phi\in [0.5, 1)$, the right side above is less than 0, we get
$$
\frac{1}{2\Delta t}(\parallel\phi^{n+1}\parallel_h^2-\parallel\phi^{n}\parallel_h^2)\leq0.
$$
Hence, inequality (\ref{eq15}) holds and the proof is complete.\\

\section{Numerical simulation}
In this section, we will conduct some numerical experiments. As the nonlinear scheme (14) which requires Newton's iteration in each step is computationally a very demanding task, we exploit the second-order semi-implicit linear scheme to  solve equation (8) with difference scheme given in (9). The noise term is just introduced in initial time for carrying the system out of metastable state. After we have done simulation for some time, the noise term will not be used, i.e., we employ the semi-implicit scheme (12) to solve equation (9) in order to simulate phase separation of the system. The mesh used is square grids of $N_x=64$ and $N_y=64$. We set $M_0=0.2$, $\chi=0.4$, $\tau =10^7$, $N=800$ and $\rho=1.0$ as in [22]. The tolerance of the conjugate gradient method is set to be $10^{-6}$. We will conduct experiments under constant time step $\Delta t=0.001$. At first, we set $T=1$, and the initial values $\phi$ are random numbers around 0.65. we get the phase evolution in Figure 4.1.
\begin{figure}[H]
\centering
\begin{minipage}{0.4\linewidth}
\includegraphics[width=1\textwidth]{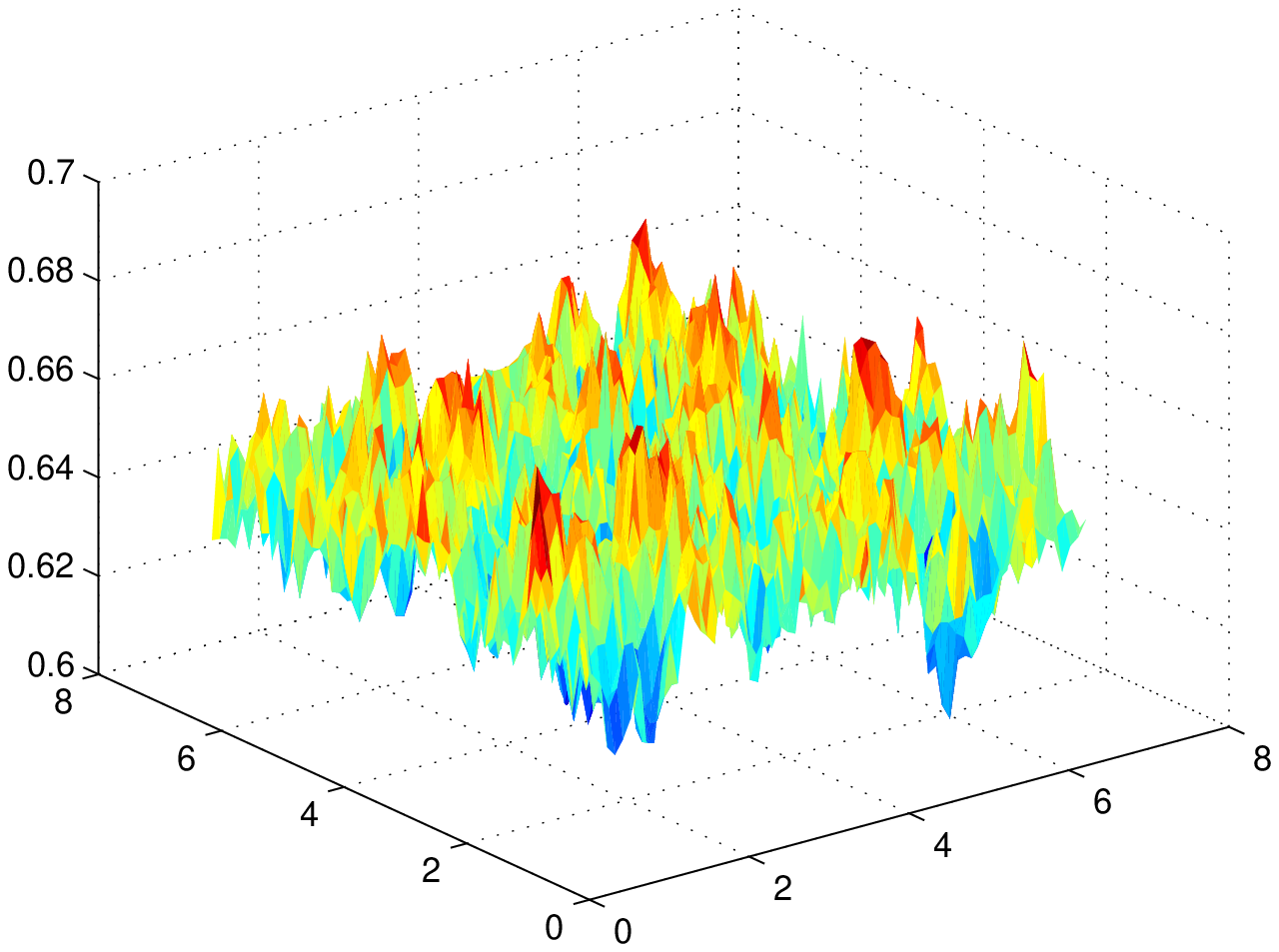}
\end{minipage}
\begin{minipage}{0.4\linewidth}
\includegraphics[width=1\textwidth]{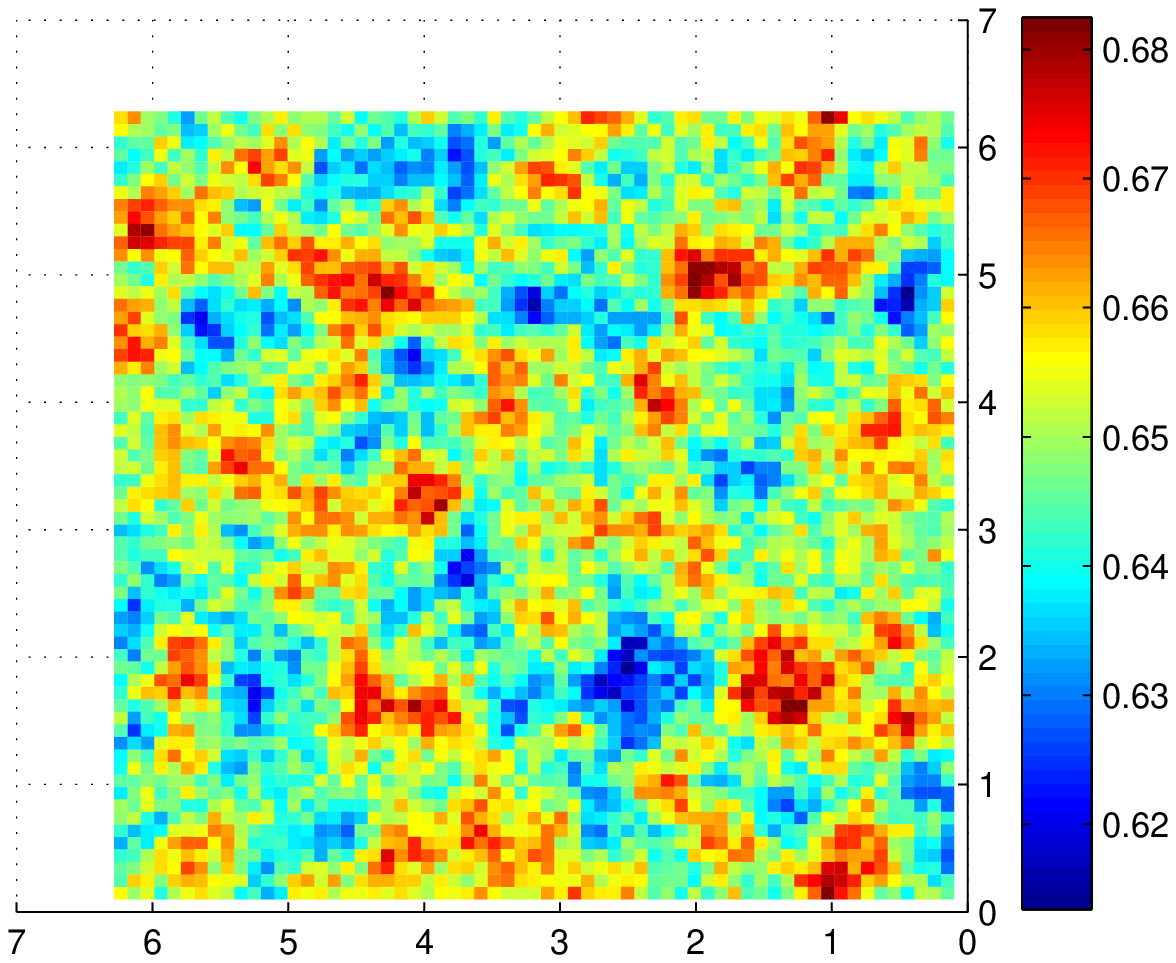}
\end{minipage}
\begin{minipage}{0.4\linewidth}
\includegraphics[width=1\textwidth]{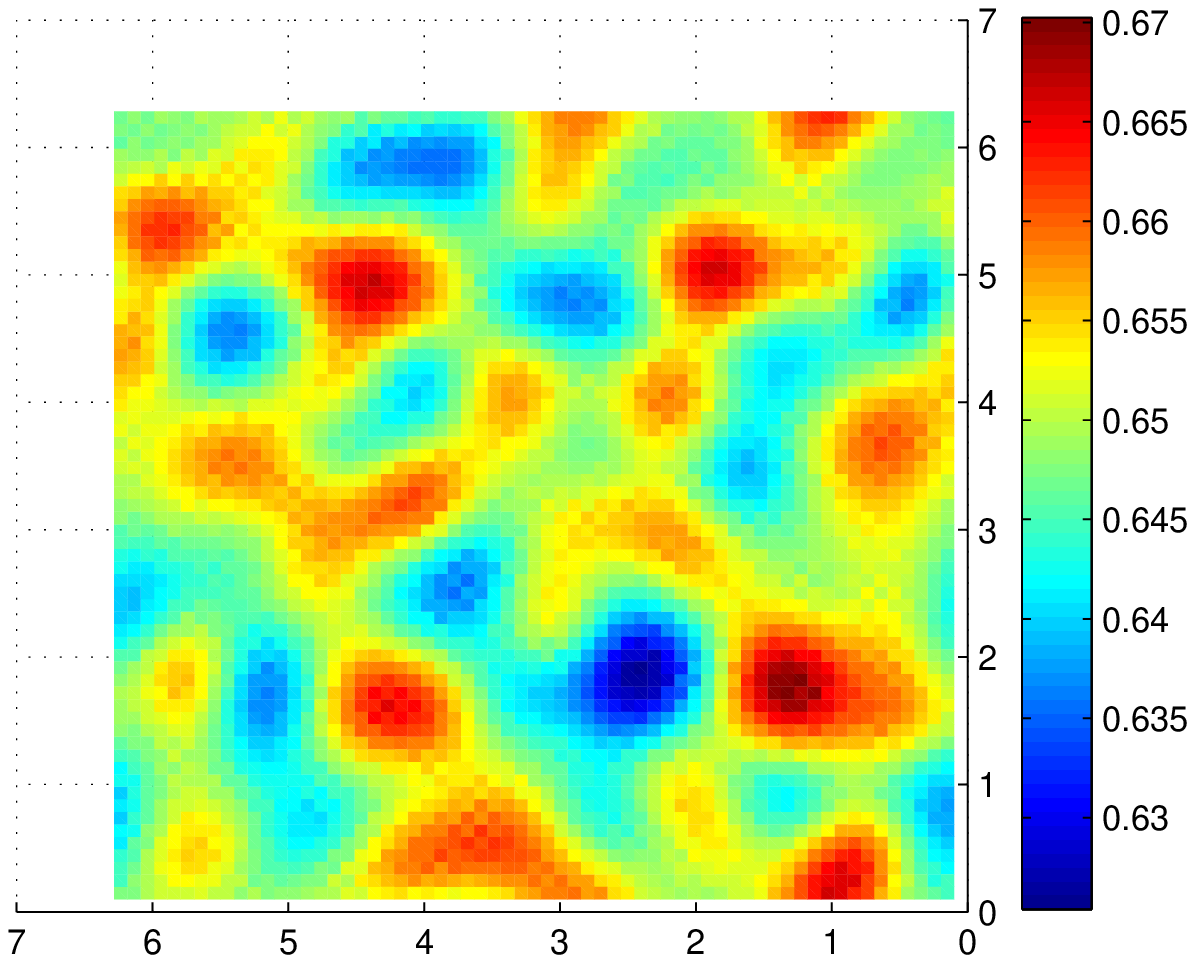}
\end{minipage}
\begin{minipage}{0.4\linewidth}
\includegraphics[width=1\textwidth]{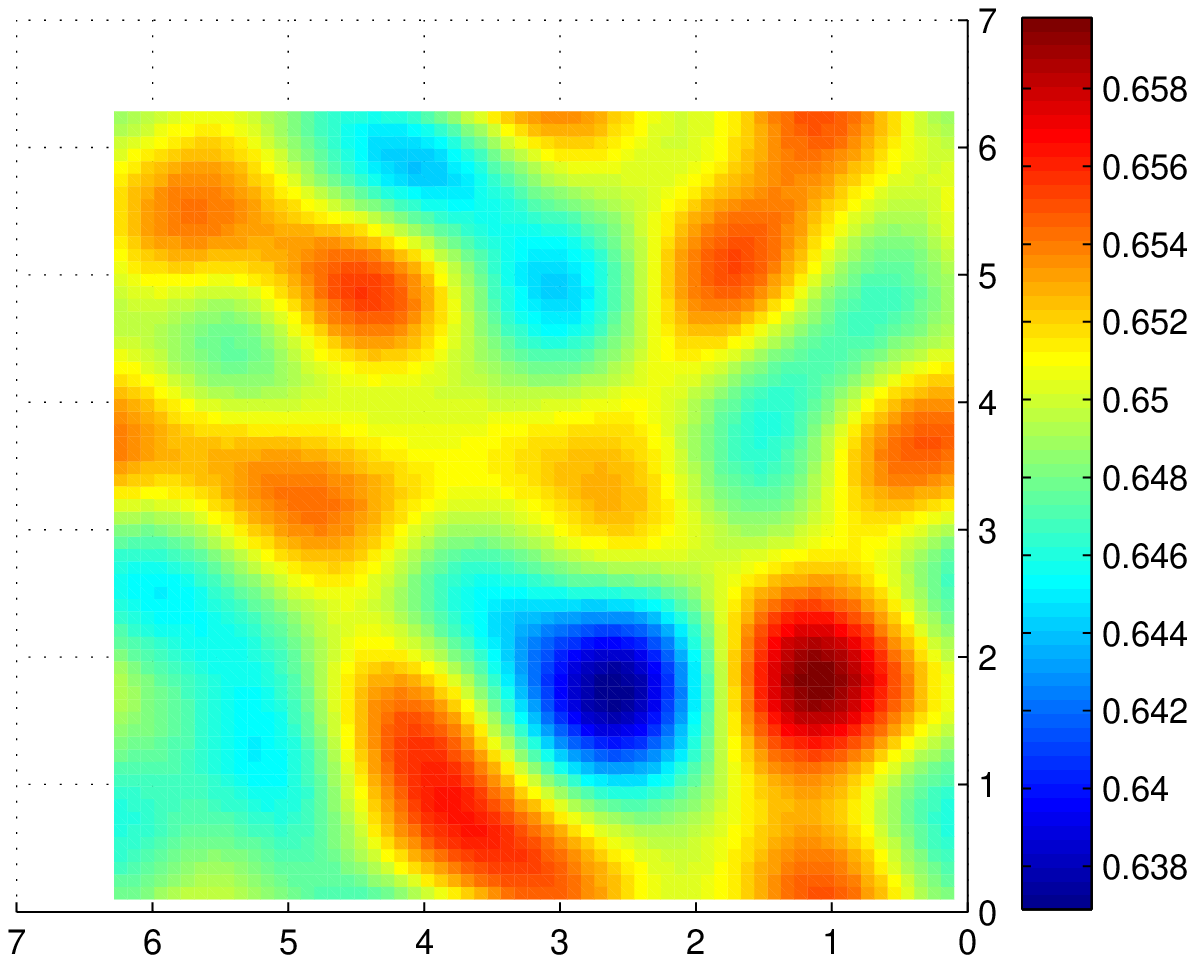}
\end{minipage}
\begin{minipage}{0.4\linewidth}
\includegraphics[width=1\textwidth]{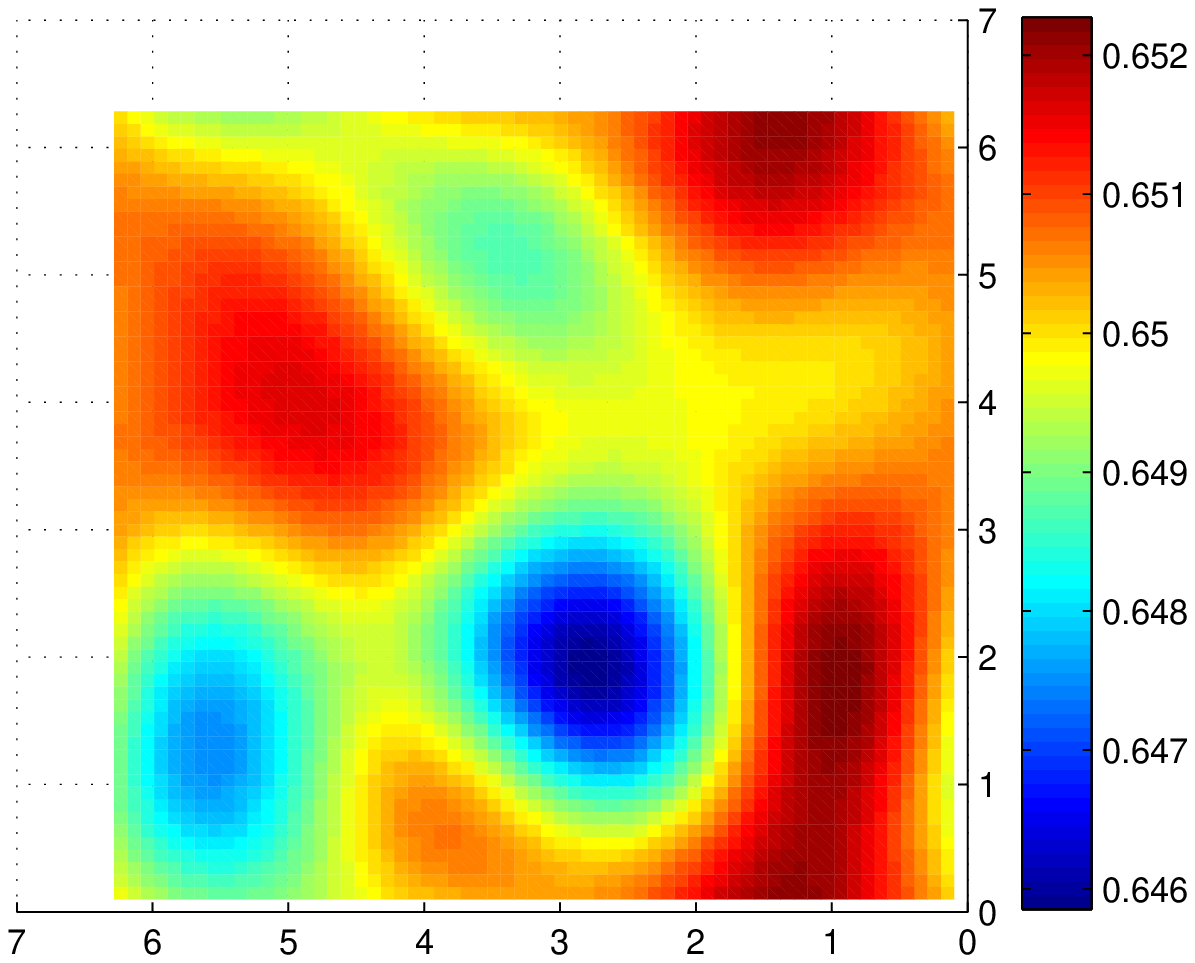}
\end{minipage}
\begin{minipage}{0.4\linewidth}
\includegraphics[width=1\textwidth]{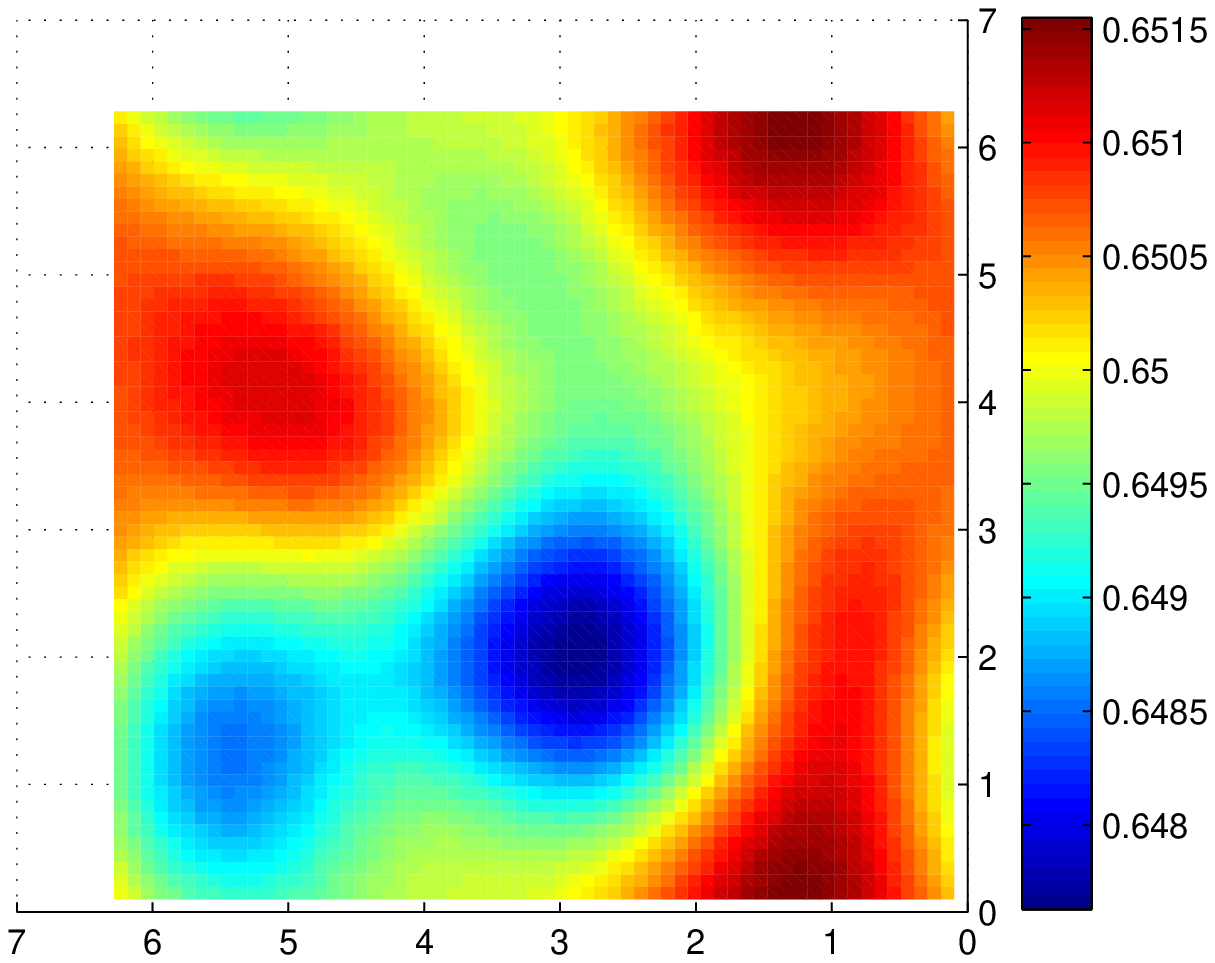}
\end{minipage}
\caption[small]{phase evolution figures at t=0.1, 10, 50, 100, 500, 1000.}
\end{figure}
From Figure 4.1, we can see phase evolution and obviously observe the phase separation. At first, hydrogel is distributed out of order in the system. As the time evolves, hydrogel gradually has well-defined structure, which will greatly improve its mechanical strength.

In the following, we change initial value $\phi$. $\phi$ are chosen as random numbers around 0.35. We get the phase evolution in Figure 4.2.
\begin{figure}[H]
\centering
\begin{minipage}{0.4\linewidth}
\includegraphics[width=1\textwidth]{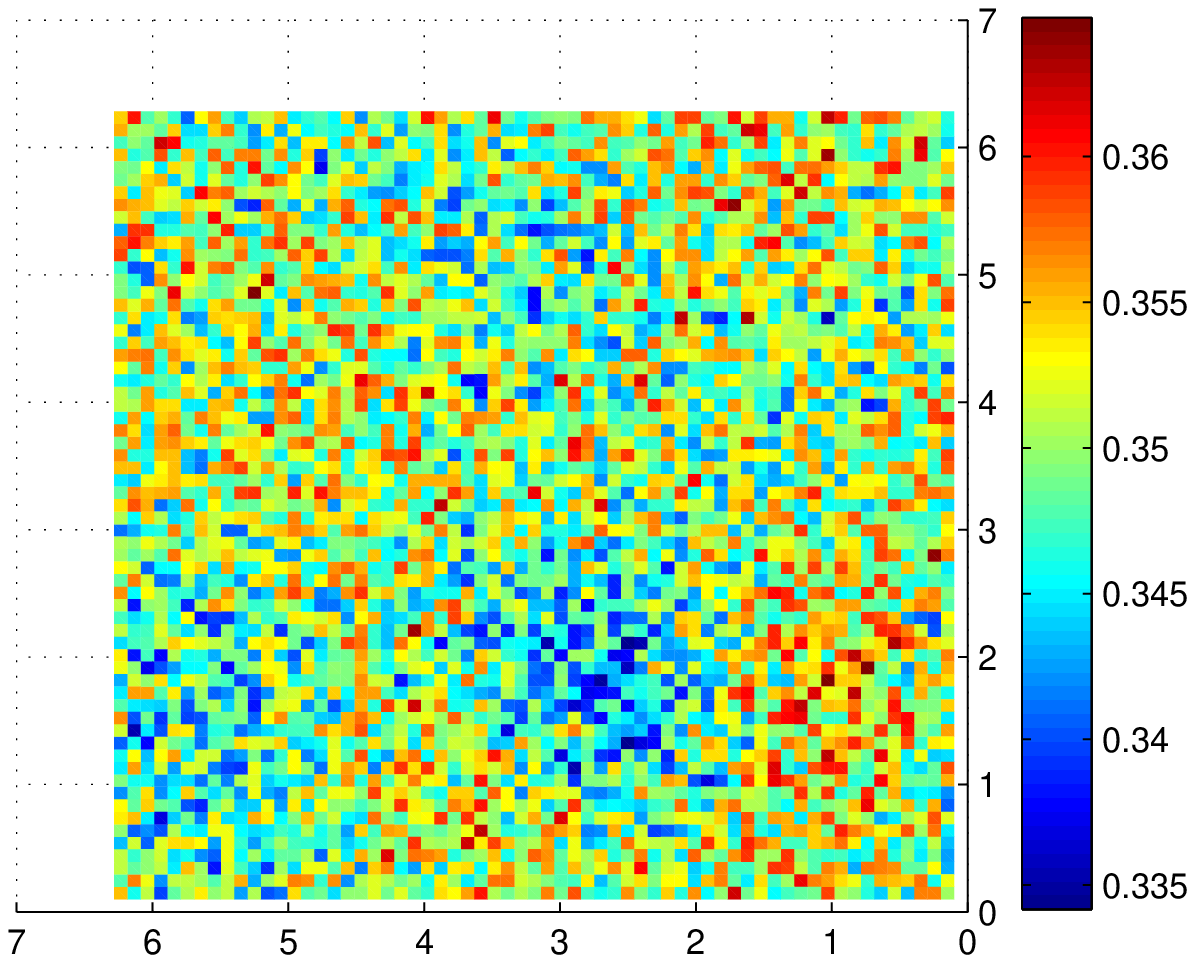}
\end{minipage}
\begin{minipage}{0.4\linewidth}
\includegraphics[width=1\textwidth]{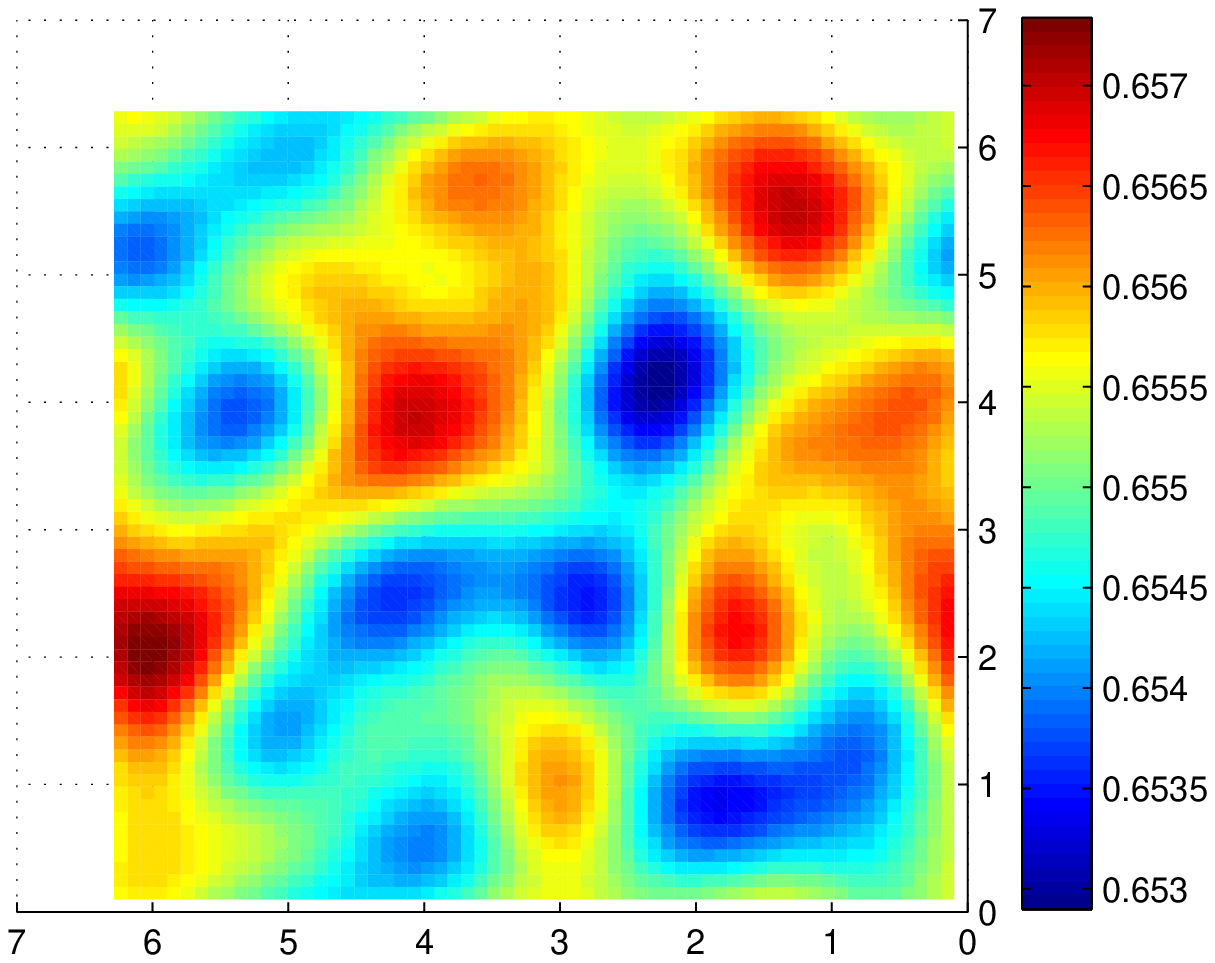}
\end{minipage}
\begin{minipage}{0.4\linewidth}
\includegraphics[width=1\textwidth]{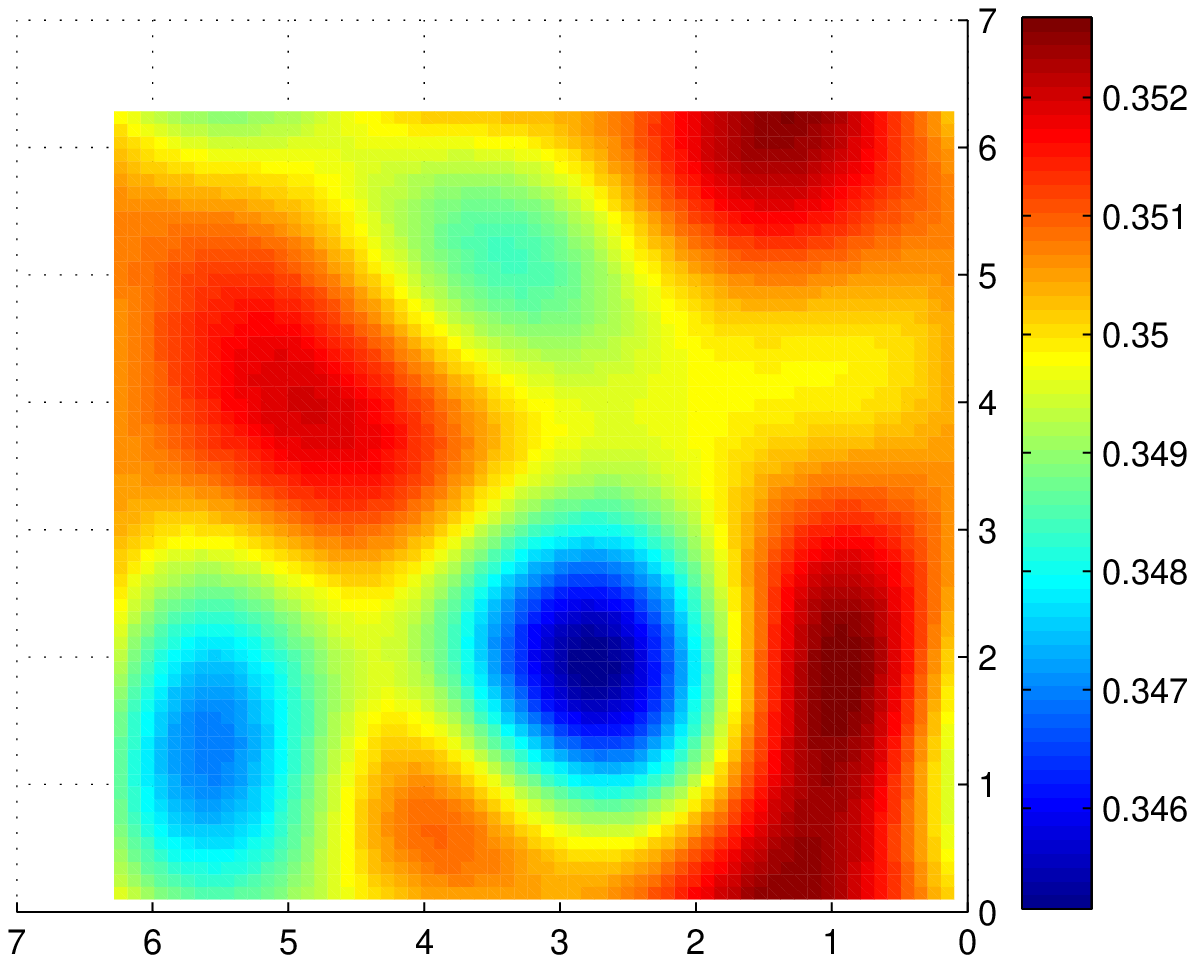}
\end{minipage}
\begin{minipage}{0.4\linewidth}
\includegraphics[width=1\textwidth]{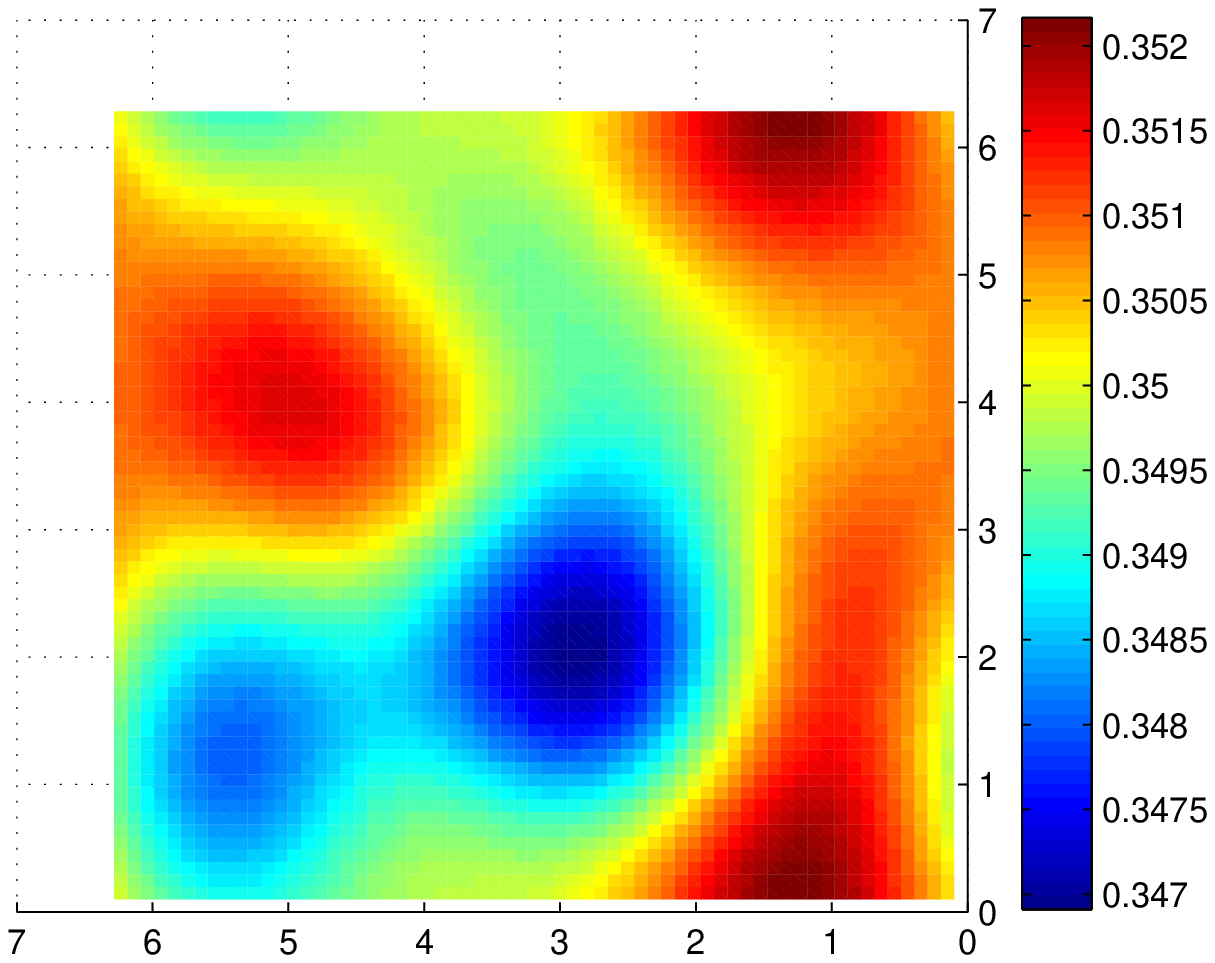}
\end{minipage}
\begin{minipage}{0.4\linewidth}
\includegraphics[width=1\textwidth]{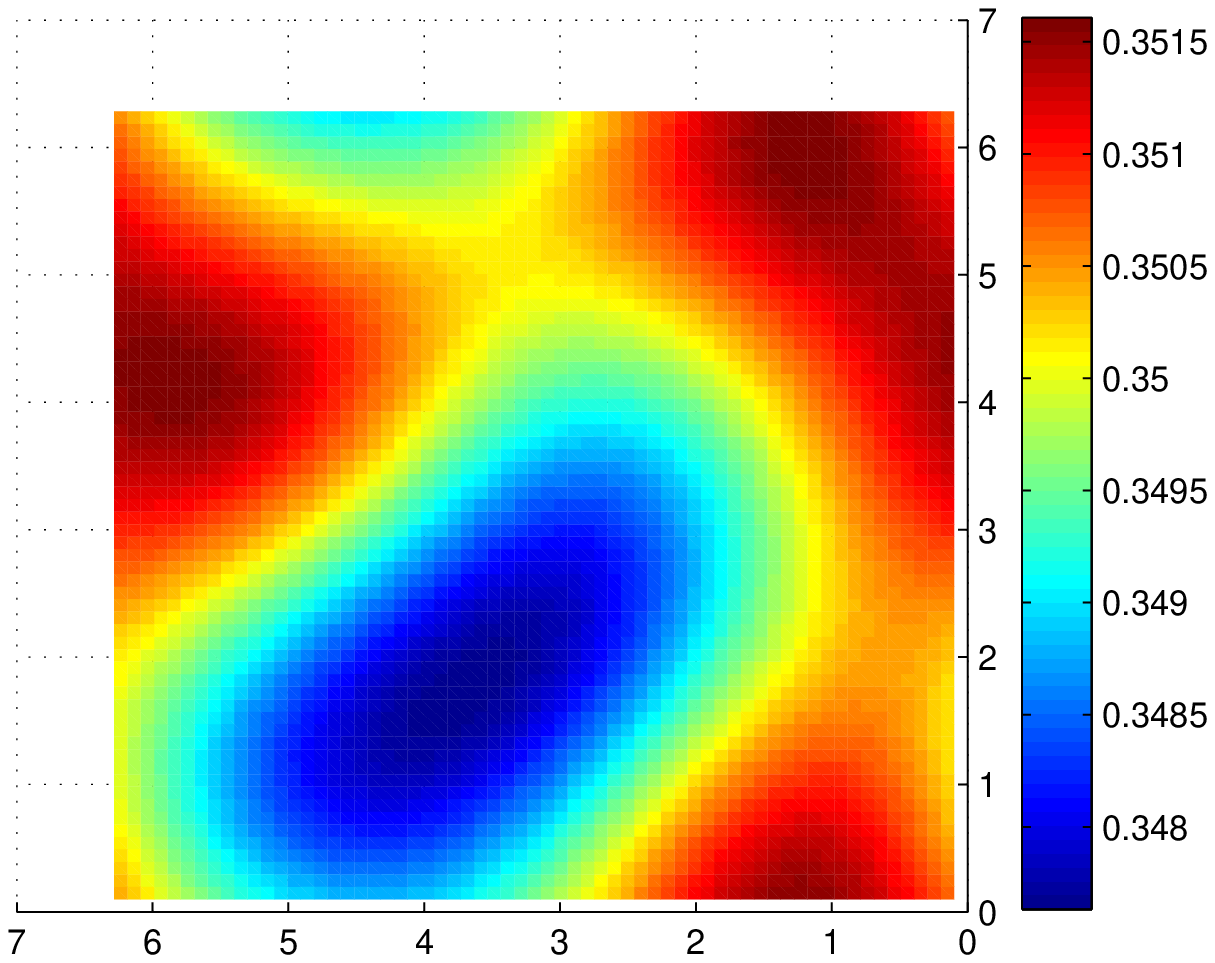}
\end{minipage}
\begin{minipage}{0.4\linewidth}
\includegraphics[width=1\textwidth]{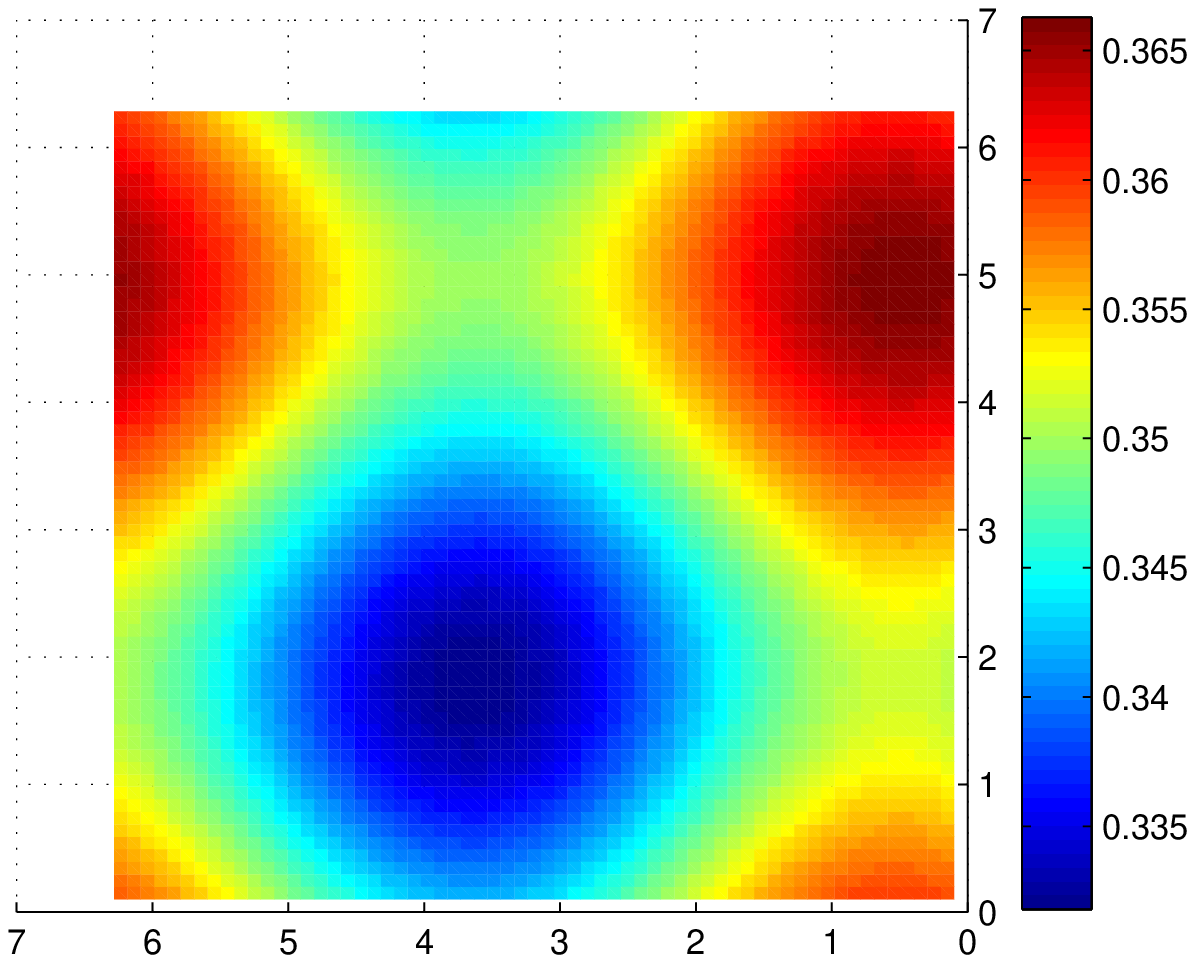}
\end{minipage}
\caption[small]{phase evolution figures at t=0.1, 10, 100, 500, 1000, 2000.}
\end{figure}
In Figure 4.2, it is showed that the phase separation of hydrogel is observed. At first, the initial value $\phi$ is considerably disordered, indicated in first phase picture in Figure 4.2. In the end, it can still evolve to a steady state with regular structure.

Meanwhile, we will change the temperature constant parameter to do some numerical experiment. The initial values $\phi$ are the same as initial values given in Figure 4.1. We just increase temperature constant to $T=20$ and T=50 respectively. The values of other parameters are same as values given above. We have the steady phase state in Figure 4.3.
\begin{figure}[H]
\centering
\begin{minipage}{0.4\linewidth}
\includegraphics[width=1\textwidth]{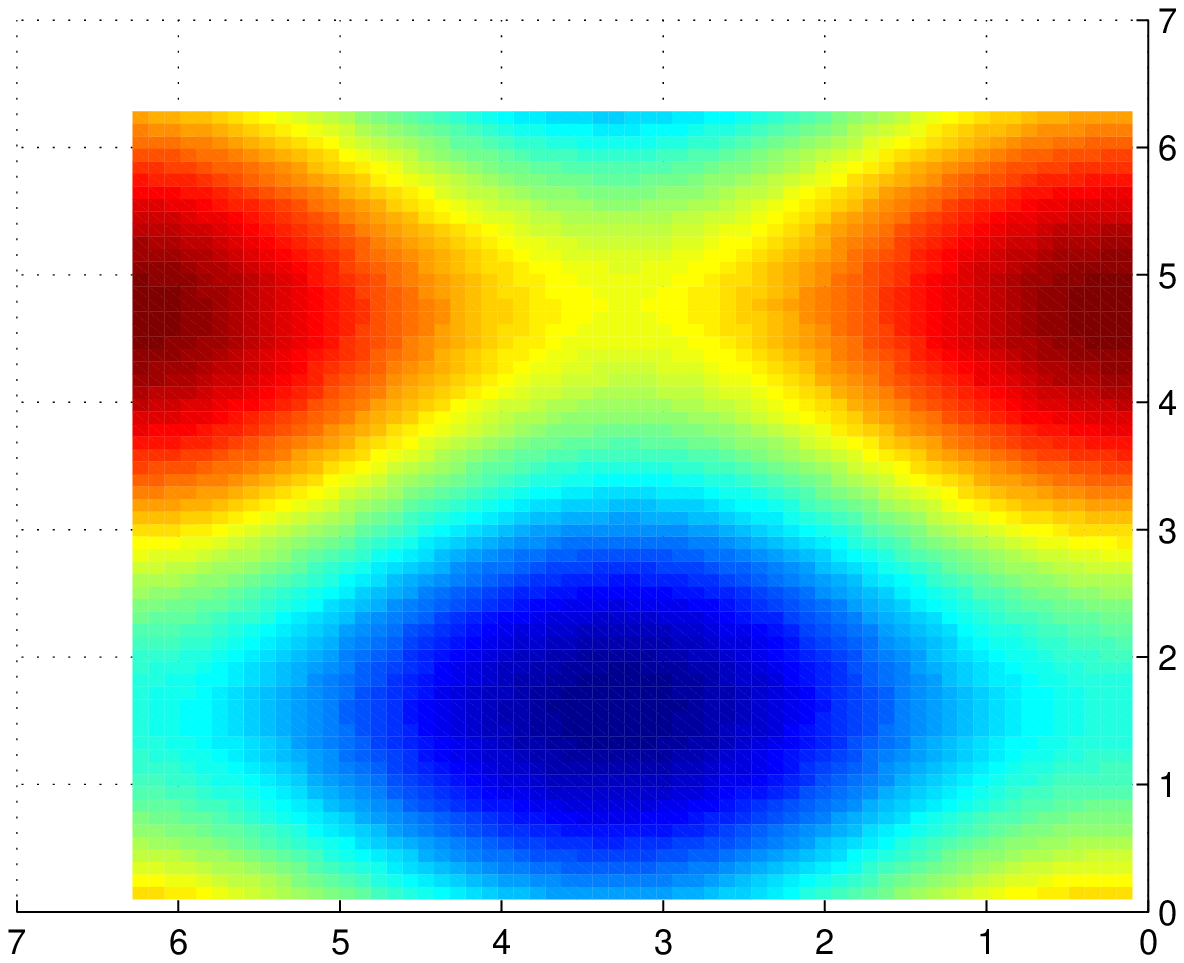}
\end{minipage}
\begin{minipage}{0.4\linewidth}
\includegraphics[width=1\textwidth]{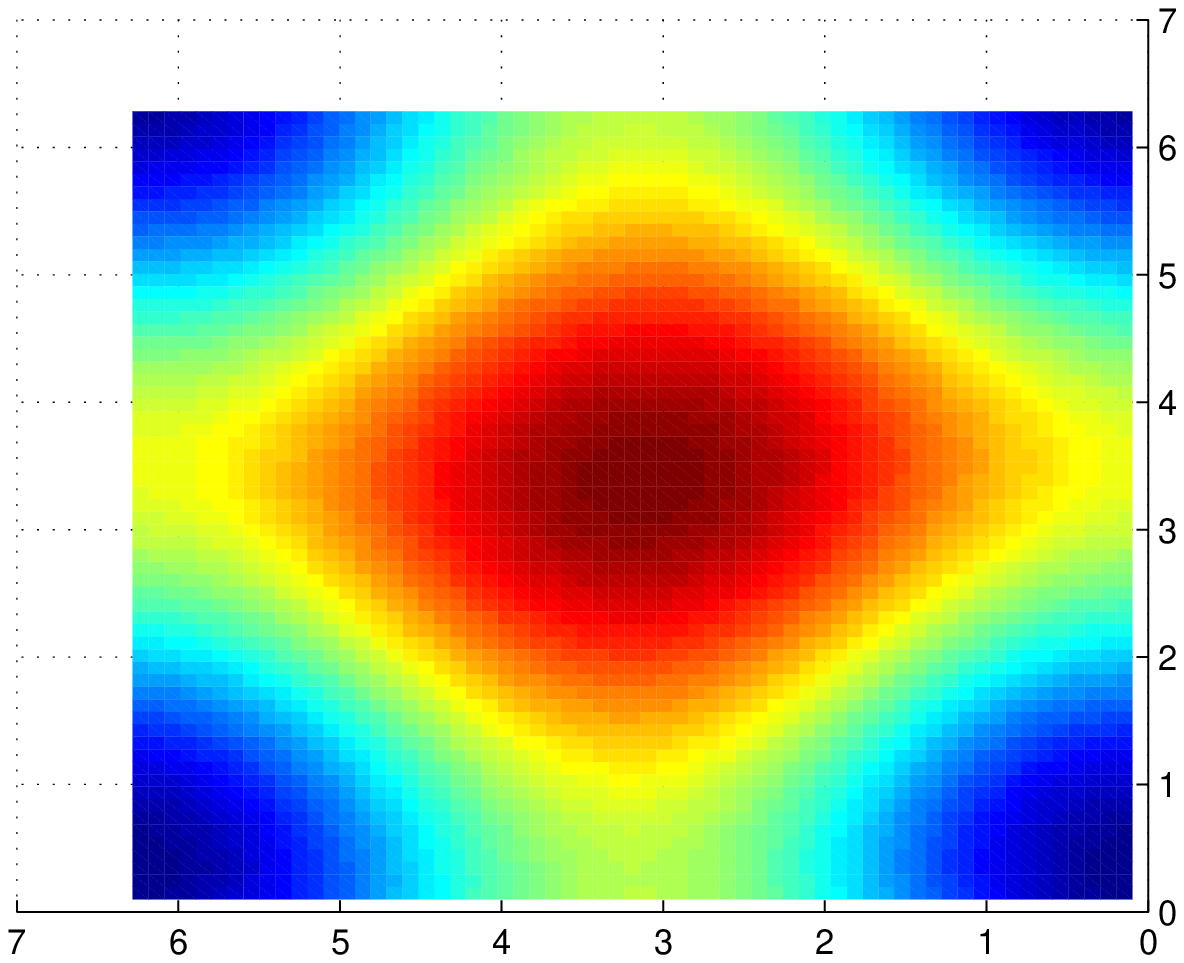}
\end{minipage}
\caption[small]{t=1000, phase figures when $T=20, 50$.}
\end{figure}
Comparing Figure 4.1 to figure 4.3, we can find that even though we have the same initial values for $\phi$, different temperature constants can derive different steady phase figures. All of these experiments where clear phase separation has been observed demonstrate that the scheme (\ref{12}) constructed is reasonable. It can simulate the process of phase evolution over long time simulation.

But during each step of numerical simulation, we should solve a linear system of equations in very high dimensions, for example, in our simulation the matrix is 4096 orders. Although the matrix is a sparse diagonally-dominant matrix, which can be solved by the conjugate gradient method in order to reduce the computation, it still costs a large amount of time when we will conduct a long time simulation. Especially, as illustrated in [22], in industrial applications, we have to enlarge the computational area. If we enlarge the area to $(0, 4\pi)\times(0, 4\pi)$, the dimensional of matrix is 4 times than that above when we have to employ the same mesh size to ensure accuracy. Therefore, the computation is unacceptable for long time simulation. This stipulates us to seek method to decrease computational cost. During the simulation, we find that the total free energy has a significant variation in initial time and has a micro change in late time, which means the solution also has such property. Hence, it allows us to introduce adaptive time step as the scheme (\ref{eq12}) keeps the property of energy decreasing and is $L^2$ stable [3].\\

\section{An adaptive time stepping method}
In section 2.1, we have demonstrated that the whole system is unconditionally energy stale. And The semi-implicit linear scheme (\ref{eq12}) and nonlinear scheme are $L^2$ stable. These properties of stability allow us to introduce large time step during numerical simulations [16, 17, 18]. For the sake of accuracy, a very large time step is unacceptable except in the time intervals where solutions change considerably little. In the following simulations, we will conduct numerical simulations under different constant time steps, such as $\Delta t =0.001$, $\Delta t =0.005$, $\Delta t =0.01$ and $\Delta t =0.05$. In the following simulations, the initial values are the same values as given in figure 1. In order to test the numerical accuracy, we take the numerical solution obtained using $\Delta t=10^{-4}$ as "exact" solution.  We define $L^2$ relative error (RE) as
\begin{equation}\label{eq16}
RE=(\sum_{i=1}^{N_x}\sum_{j=1}^{N_y}(\phi^{\Delta t}_t(i,j)-\phi^{\Delta t_0}_t(i,j))^2h_xh_y)^{\frac{1}{2}},
\end{equation}
where $\phi^{\Delta t_0}_t(i,j)$ is the "exact" solution $\phi$ we assume at grid $(i,j)$ under $\Delta t=0.001$, and similarly, $\phi^{\Delta t}_t(i,j)$ is the solution at $\Delta t$($\Delta t\geq 0.001$). Here $h_x=\frac{2\pi}{N_x}$, $h_y=\frac{2\pi}{N_y}.$
We set constant temperature T=1. The following figures show total free energy evolution over time:
\begin{figure}[H]
\centering
\begin{minipage}{0.4\linewidth}
\includegraphics[width=1.1\textwidth]{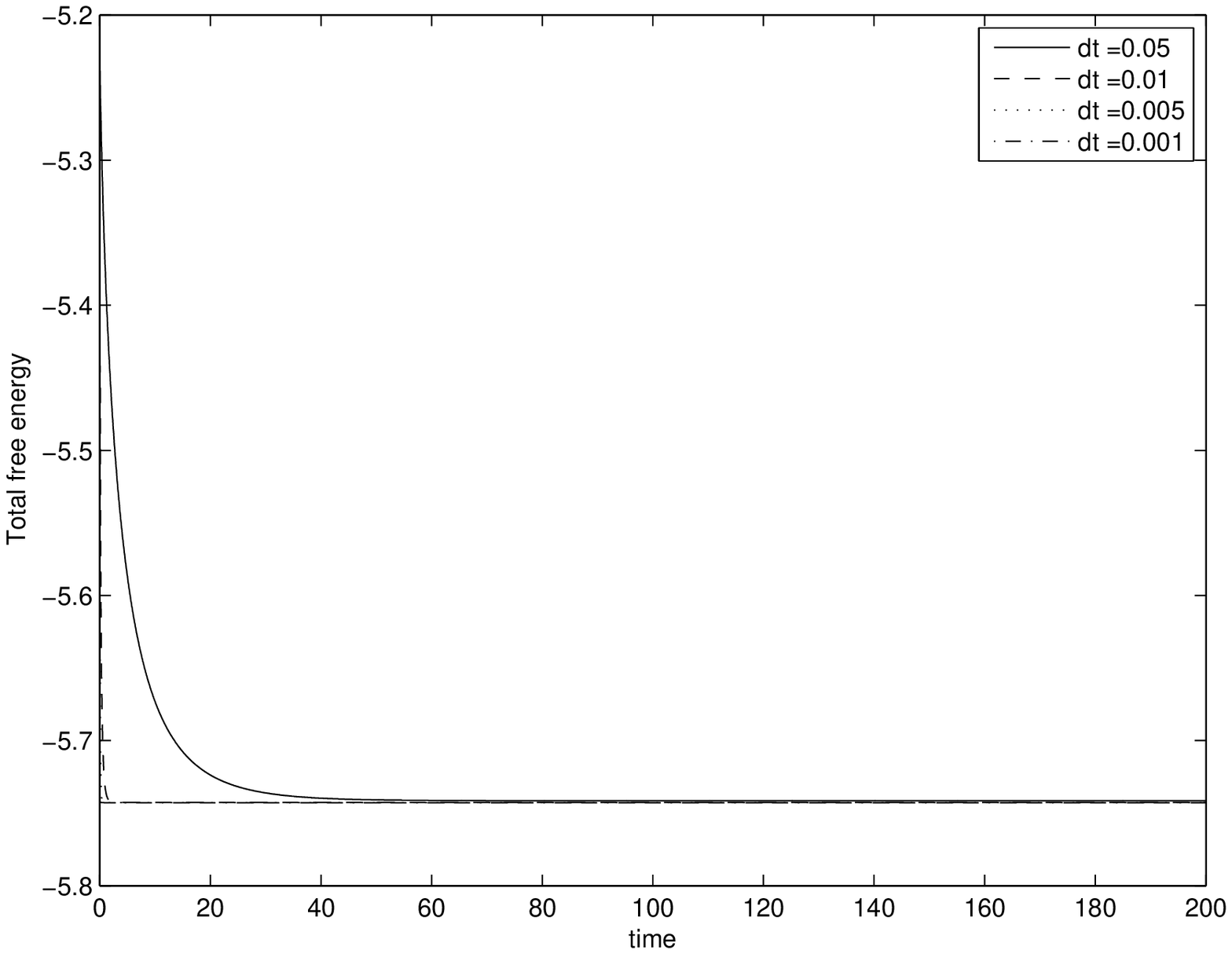}
\end{minipage}
\begin{minipage}{0.4\linewidth}
\includegraphics[width=1.1\textwidth]{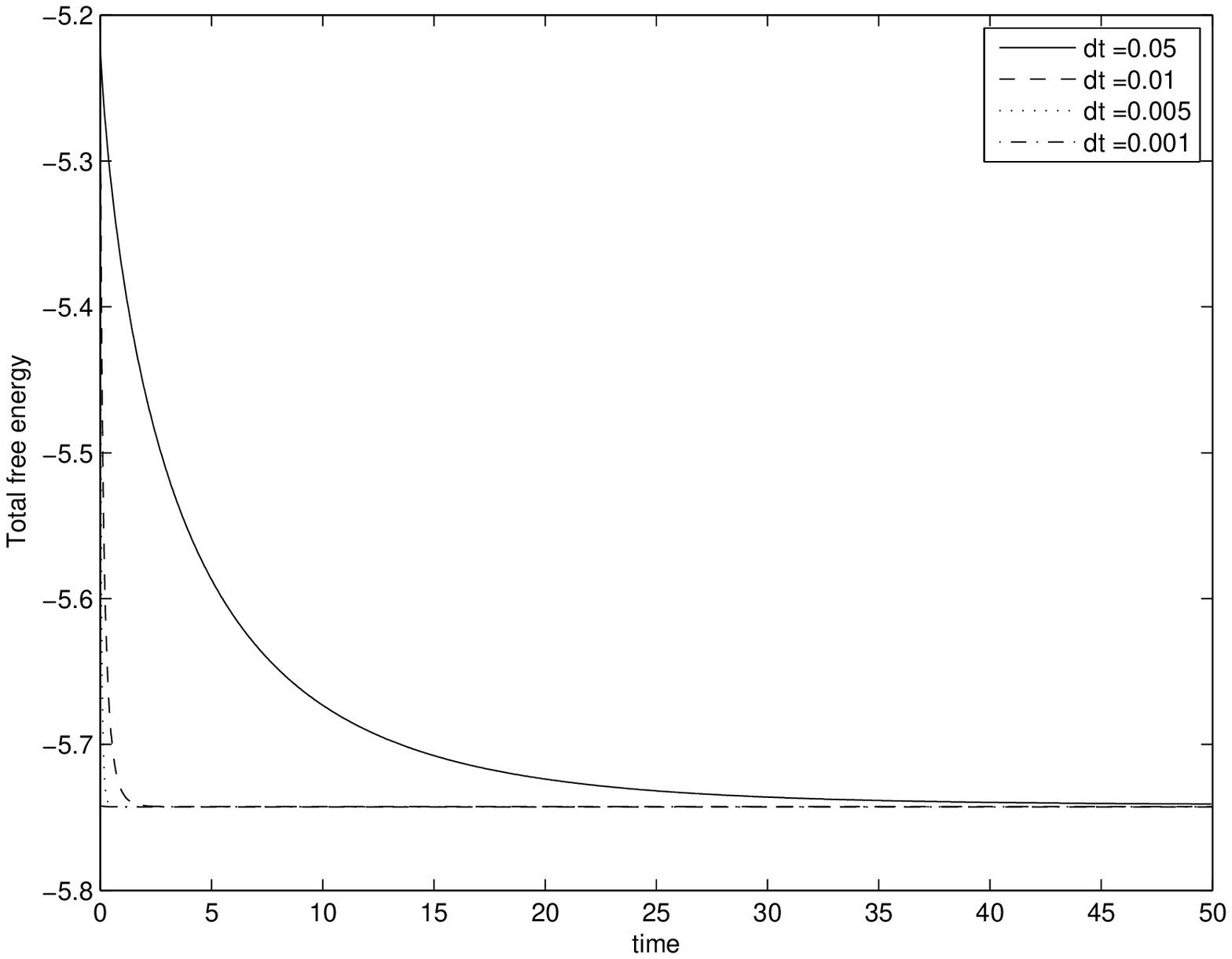}
\end{minipage}
\begin{minipage}{0.4\linewidth}
\includegraphics[width=1.1\textwidth]{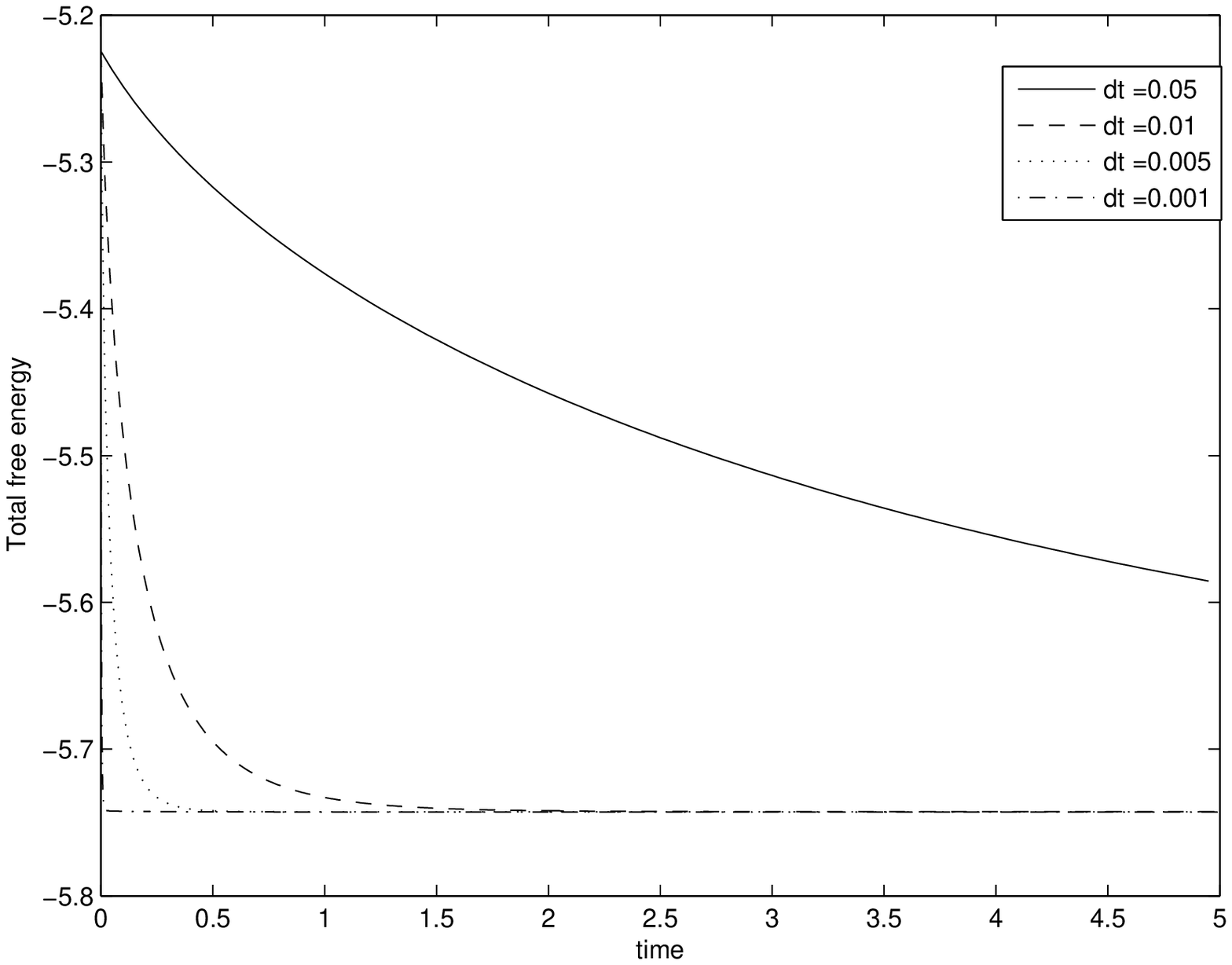}
\end{minipage}
\begin{minipage}{0.4\linewidth}
\includegraphics[width=1.1\textwidth]{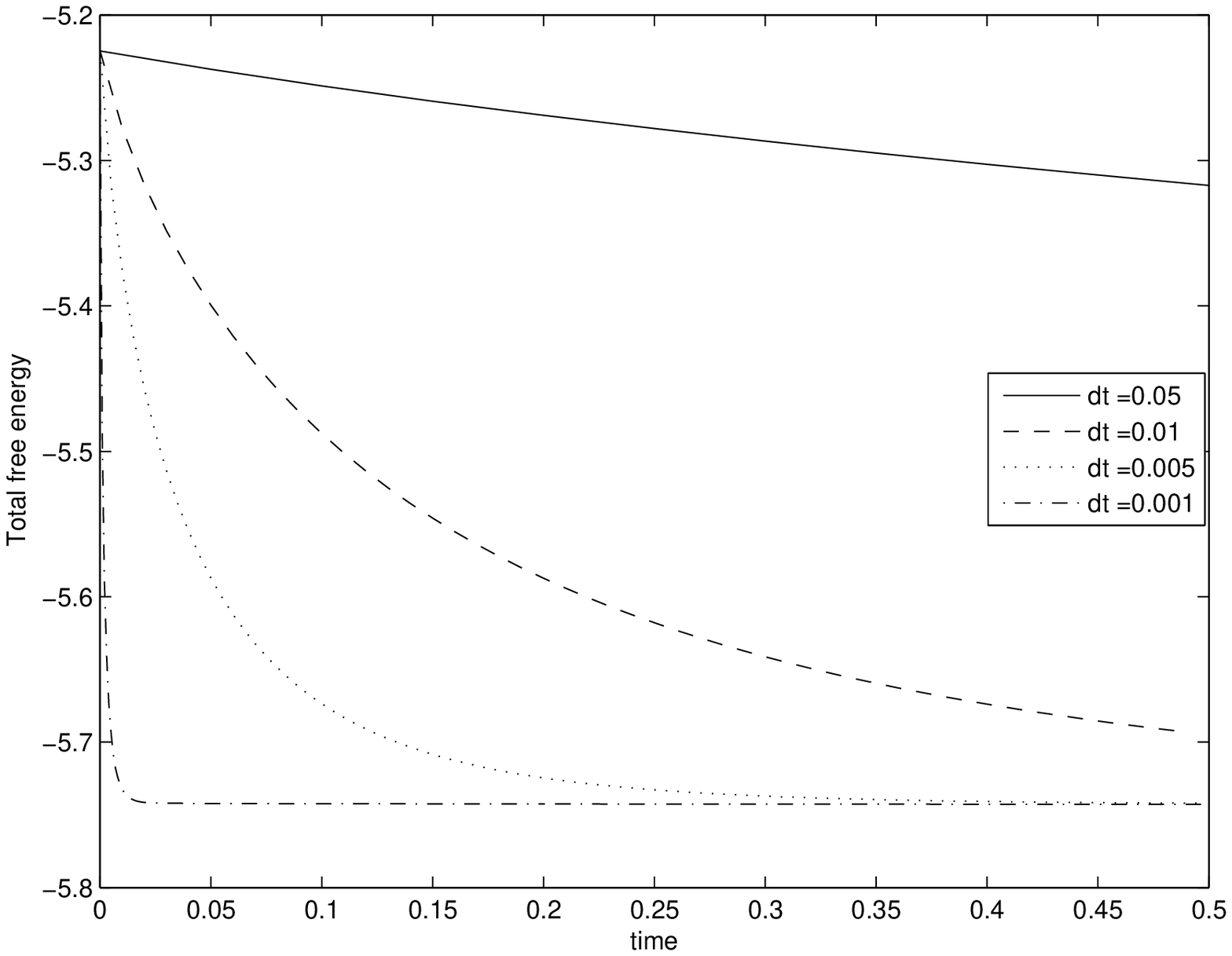}
\end{minipage}
\caption[small]{Total energy evolution over time at T=1 under different constant time step.}
\end{figure}

From Figure 5.4, we can see total free energy under different constant time steps are all decreasing over time. This indicates that numerical solutions under different constant time steps are consistent with the property of energy stability analyzed in section 2.1. Although there are errors at first, as time evolves, solutions under different constant time steps are the same ultimately. In initial time, solutions under larger constant time steps have higher errors than solutions under small constant time steps. This demonstrates that the systems are sensitive to time steps chosen at initial time. However, we can also see from figures that the total free energy decreases significantly at initial time and changes minor after some time. As we have proven that the whole system is energy stable and the scheme \ref{eq12} is $L^2$ stable, this allows us to choose large time step without resulting in much error.

At the same time we will conduct some numerical experiments under different temperature to test the robustness of observations above. In the following numerical simulations, we choose $T=50$.
\begin{figure}[H]
\centering
\begin{minipage}{0.4\linewidth}
\includegraphics[width=1.1\textwidth]{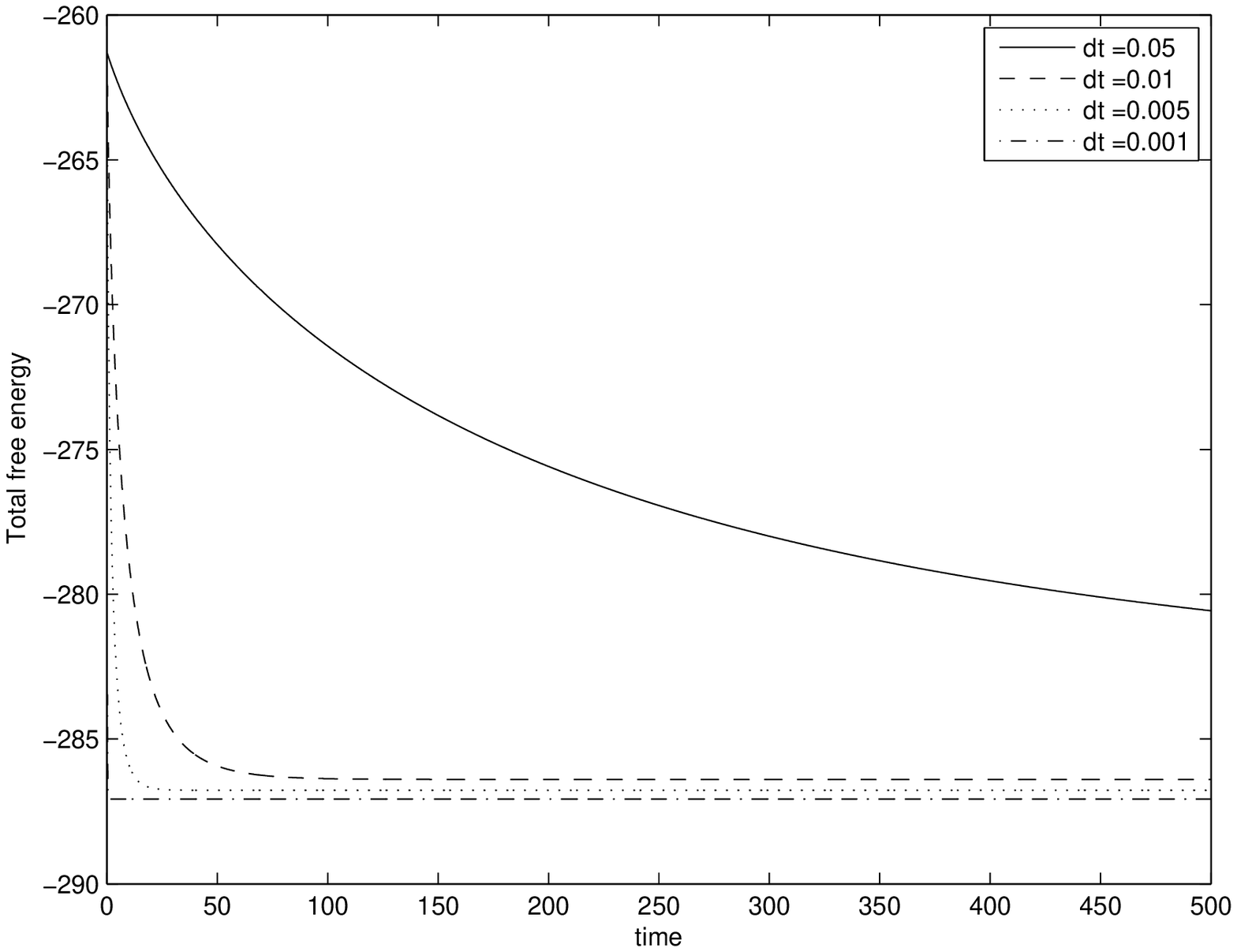}
\end{minipage}
\begin{minipage}{0.4\linewidth}
\includegraphics[width=1.1\textwidth]{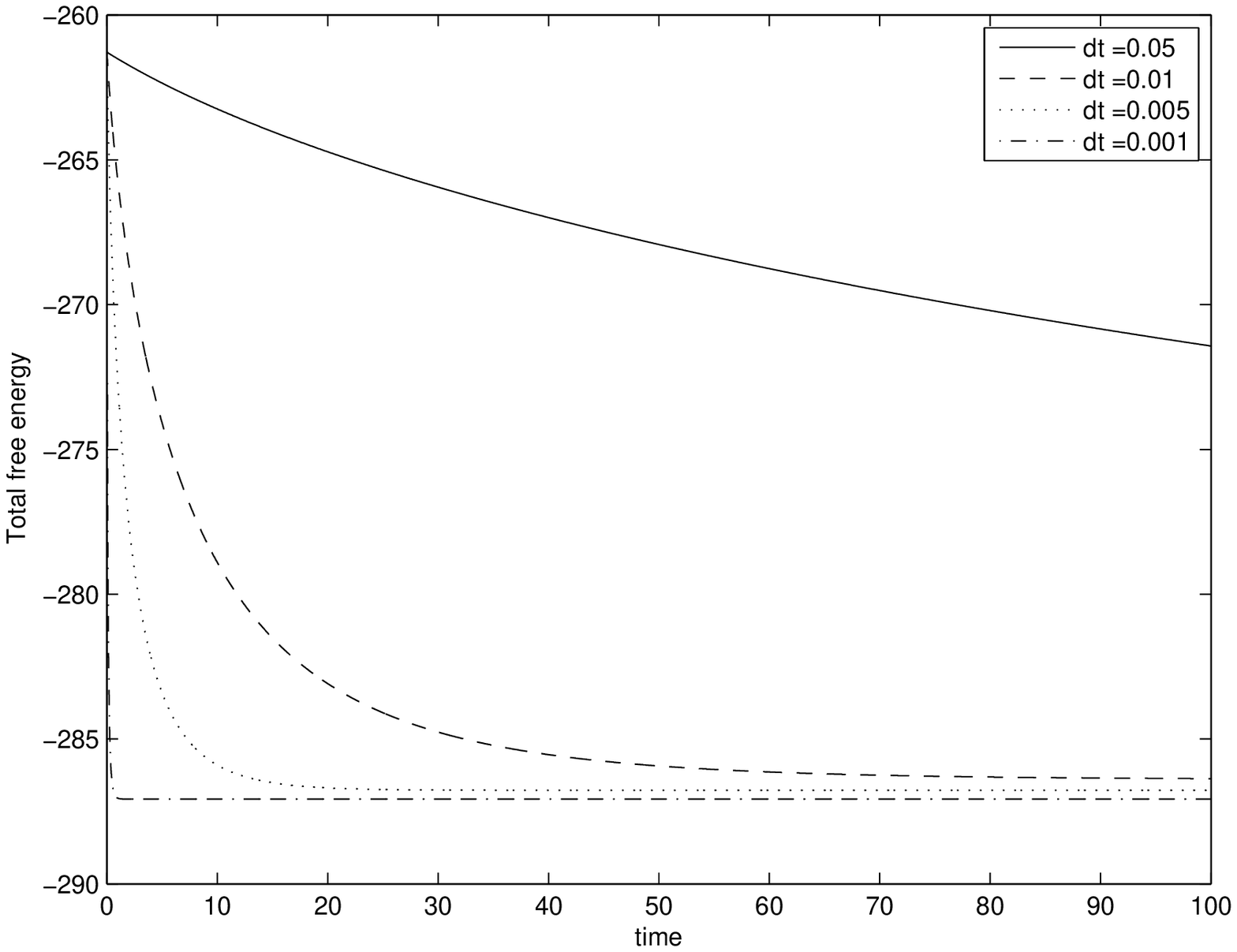}
\end{minipage}
\begin{minipage}{0.4\linewidth}
\includegraphics[width=1.1\textwidth]{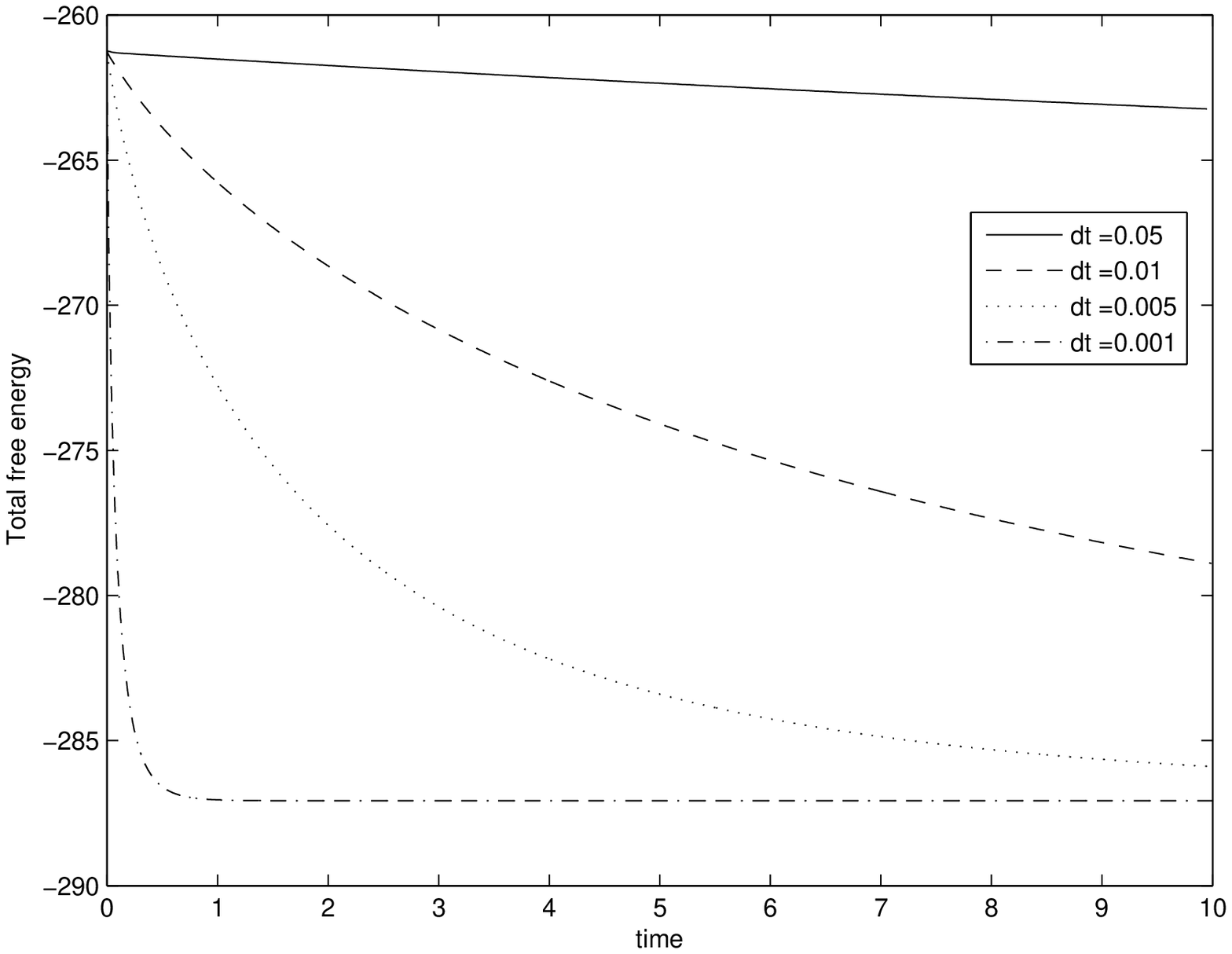}
\end{minipage}
\begin{minipage}{0.4\linewidth}
\includegraphics[width=1.1\textwidth]{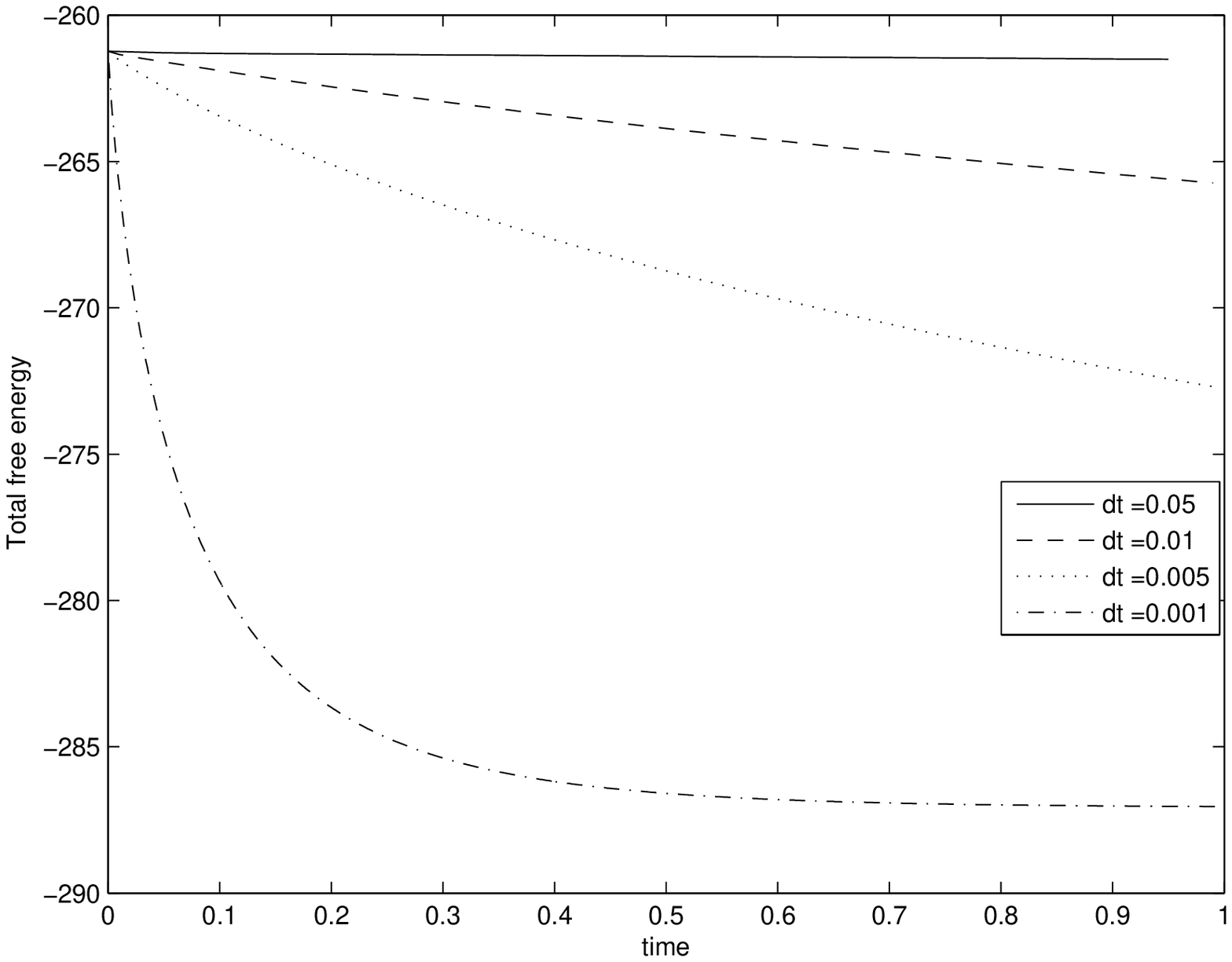}
\end{minipage}
\caption[small]{Total energy evolution over time at T=1 under different constant time step.}
\end{figure}

From Figure 5.5, we can see that when temperature increases, large constant time step will deviate the energy curve from the energy curve obtained by "exact" solution, which illustrate that the solution obtained by larger constant time step is much higher than that of solution with small temperature under the same initial value condition . This indicates that numerical solutions are more sensitive to constant time steps. Therefore, we should be more concerned about time steps when we numerically solve the equation.

In the following, we will introduce an adaptive time step to numerically solve equations. The basic idea of choosing adaptive time step is that when the total free energy decreases significantly, which mean solutions have a large variation, we employ small time step and when total free energy has a micro change, we will choose large time steps for it will not bring about large error under this circumstance. In this way, small errors can be guaranteed during the whole process.

Following [3], according to the choice of adaptive time step and definition of total free energy in this equation, we improve the choice of adaptive time step at time t as:
\begin{equation*}
\Delta t=\max(\Delta t_{min}, \lambda(T,t)\frac{\Delta t_{max}}{\sqrt{1+\mu |U'(t)|^2}}).
\end{equation*}
Here $\mu$ is a constant parameter, $U'(t)$ is the derivative of $U(t)$ and $\lambda(T,t)$ is introduced to ensure accuracy and efficiency.
For the same change of total free energy, we will choose small time step at first for accuracy but it is safe to use a larger time step in late time for the system evolves to a steady state. The total free energy is proportional to temperature T. For the same change of $\phi$, the derivative of $\mu |U'(t)|^2$ at $T=50$ is $2500$ times than the derivative of it at T=1. In this case, for the same change of $\phi$, the choice of adaptive time step is considerably small. So we have to multiply $\lambda(T,t)$ to ensure the time step is not too small.

In the following simulation, we choose $\Delta t_{min}=0.001,\Delta t_{max}=0.1,\mu=1000$. and the definition of $\lambda(T,t)$ is as follows.
\begin{equation*}
\lambda(T,t)=
\begin{cases}
1,\quad if\quad 0\leq t\leq 100.\\
1.5,\quad if\quad 100< t\leq 200.\\
2,\quad if\quad 200< t\leq 300.\\
3,\quad if\quad 300<t\leq 400.\\
4,\quad if\quad 400<t\leq 500.\\
5,\quad if\quad t> 500
\end{cases}.
\end{equation*}
First, we list the CPU time consumed at some selected time levels $t=1$, $t=10$, $t=100$, $t=1000$ respectively under different time steps.
\begin{center}
\title{CPU time costs under different time steps at $t=1$, $t=10$, $t=100$, $t=500$}.\\
\begin{tabular}{|c|c|c|c|c|}
               \hline
               Running time&t=1 &t=10 & t=100 & t=500 \\ \hline
               $\Delta t=0.05$& 47.8 &388.6 & 3416.9&14341.4 \\ \hline
               $\Delta t=0.01$&149.6& 905.6 & 4190.2 & 5987.2 \\ \hline
               $\Delta t=0.005$&151.8 & 948.5 & 2335.1 & 3920.5 \\ \hline
               $\Delta t=0.001$&184.5 & 419.7 & 2125.2 & 10485.2 \\ \hline
               The adaptive time step&110.1 & 202.7 & 668.8 & 2809.4 \\
               \hline
             \end{tabular}
\end{center}

From the table above, we can see that the CPU time consumed under adaptive time steps is much less than that of time under other constant time steps. Meanwhile, we should consider the accuracy of solutions evolved. At first, we see the
\begin{figure}[H]
\centering
\begin{minipage}{0.4\linewidth}
\includegraphics[width=1.1\textwidth]{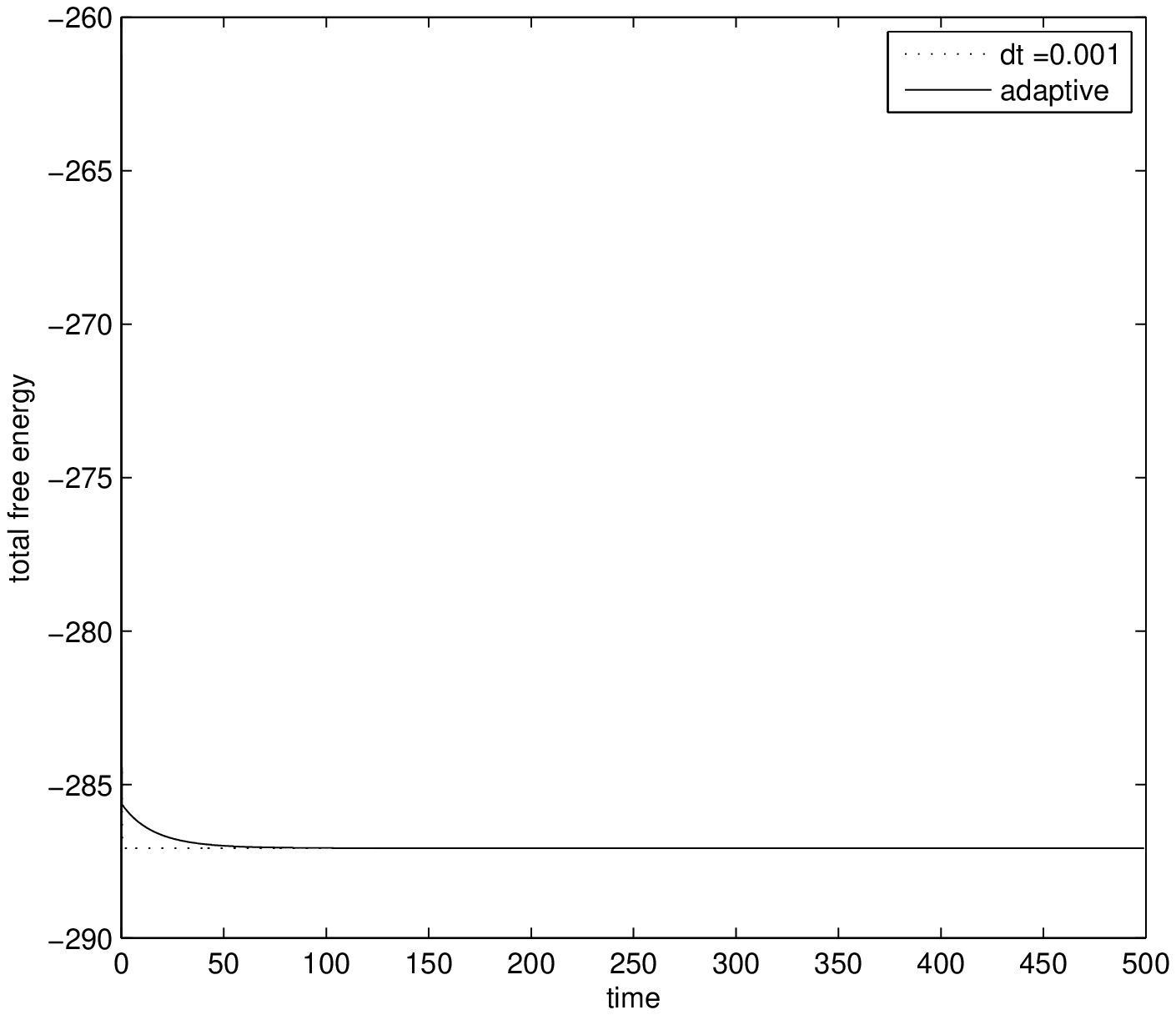}
\end{minipage}
\begin{minipage}{0.4\linewidth}
\includegraphics[width=1.1\textwidth]{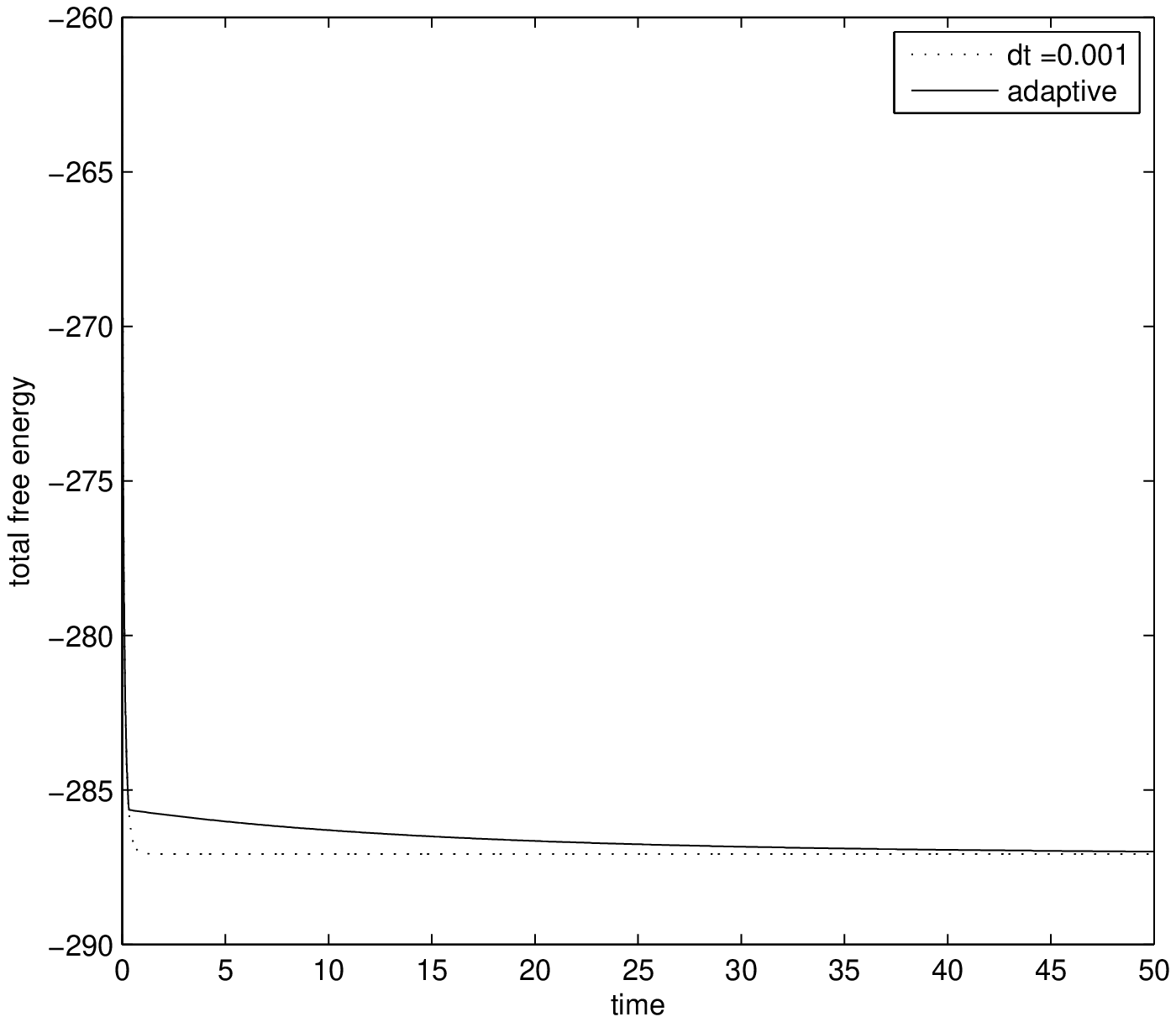}
\end{minipage}
\caption[small]{Total free energy evolution over time, blue:$\Delta t=0.001$, red: the adaptive time step.}
\end{figure}

From Figure 5.6, we can see that the difference of total free energy evolution is very small under constant time step $\Delta t=0.001$ and
the adaptive time step. The main reasons are as follows. At the initial time with significant change of total free energy, the choice of time step in the adaptive method is almost the same with that of time step in constant time step simulation $\Delta t=0.001$. However, as the system evolves for some time, the free energy changes little indicating that the change of $\phi$ is very small. Under this circumstance, we can choose a large time step which might be $100$ times than $0.001$ without bringing about much error.

Also, we will give $L^2$-errors to test the accuracy of different time steps at some selected time levels. In the following, we give comparison of $L^2$-errors of different time steps both constant time steps and the adaptive time step at time levels,$t=1$, $t=10$, $t=100$, $t=1000$ respectively.
\begin{center}
\title{$L^2$ error under different constant time steps and the adaptive time step at $t=1$, $t=10$, $t=100$, $t=200$}.\\
\begin{tabular}{|c|c|c|c|c|}
  \hline
   $L^2$-error& t=1&t=10 & t=100 & t=200 \\\hline
  $\Delta t=0.05$&0.031 & 0.0286& 0.0256 & 0.0376 \\\hline
  $\Delta t=0.01$ &0.0267& 0.0152 &0.0086  & 0.0083 \\\hline
  $\Delta t=0.005$ &0.0198& 0.0058 & 0.0037 & 0.0056 \\\hline
  The adaptive time step&0.0095& 0.0042 & 0.000987 & 0.000325  \\
  \hline
\end{tabular}
\end{center}

From the $L^2$ errors given above, we can see that larger constant time step have higher errors which is consistent with analysis given above. The $L^2$-error with adaptive time strategy is considerably small during the whole process of simulation, which indicates that the adaptive method is not only efficient but also accurate.\\

\section{Concluding remarks}
In this paper, we have demonstrated the simplified version of TDGL equation is $L^2$ stable. At the same time, we put forward two second-order semi-implicit schemes satisfying energy decreasing property, which are $L^2$ stale and consistent with the analytic property of $L^2$ stability. We also have conducted some numerical experiments to test the robustness of the schemes. Some numerical results are shown above using one of schemes constructed. From the numerical results, we can clearly see the phase evolution of MMC hydrogels and the process of phase separation has been observed during the experiments, which indicates effectiveness of the numerical scheme (\ref{eq12}).

Meanwhile, during the simulation, we find that the total free energy of the system changes significantly at initial time and varies little after some time simulation, which mean that numerical solutions have a palpable variation at first and have a mirror change in late time. As the system is energy stale and the schemes are $L^2$ stable with energy decreasing property, which guarantees that large time step can be introduced at late time without bring about much error. Hence we implement the adaptive time-stepping method to simulate the equation. The numerical results show the efficiency and accuracy of the adaptive time strategy. These results indicates that simplified TDGL equation is very suitable for simulating the phase evolution of MMC hydrogels after introducing the adaptive time strategy.

Although we have gotten some reasonable results to show that scheme (\ref{eq12}) is a robust scheme to simulate the phase separation, much further work should be done to test the accuracy of the scheme. As spectral method is high resolution method to numerically solve equations with periodic boundary condition, we can exploit the spectral method to compute (\ref{eq8}) and make a comparison between solutions obtained with scheme (\ref{eq12}) and solutions with spectral method under the same initial value condition. In our further effort, we will employ the spectral method to solve the equation and investigate more about TDGL equation in its initial form.


\begin{thebibliography}{99}
\bibitem{pa1} A. Chakrabarti, R. Toral, J. D. Gunton, et al; Dynamics of phase separation in a binary polymer blend of critical composition, J. Chem. Phys 1990. 92, 6899.
\bibitem{pa2} M. Fermeglia, S. Pricl; Multiscale modeling for polymer systems of industrial interest, Progress in organic coatings, 2007, 58: 187-199
\bibitem{pa3} H.W. Gibson, M.C. Bheda, P.T. Engen; Rotaxanes, catenanes, polyrotaxanes,polycatenanes and related materials. Prog Polym Sci 1994;19:843每945.
\bibitem{pa4} J.P. Gong, Y. Katsuyama, T. Kurokawa, et al; Double network hydrogels with extremely high mechanical strength. Adv Mater 2003;15:1155每8.
\bibitem{pa5} W.W. Graessley, D.S. Pearson; Stress每strain behavior in polymer networks containing nonlocalized junctions. J Chem Phys 1977;66:3363每70.
\bibitem{pa6} M. Huang, H. Furukawa, Y. Tanaka, et al; Importance of entanglements between first and second components in high-strength double network gels. Macromolecules 2007;40:6658每64.
\bibitem{pa7} K. Harag, T. Takehisa; Composite hydrogels: a unique organic每inorganic network structure with extraordinary mechanical, optical, and swelling/deswelling properties. Advance Matter.2002;14:1120-4.
\bibitem{pa8} T. Huang, H.G. Xu, K.X. Jiao, et al; A novel hydrogel with high mechanical strength: A macromolecular microsphere composite hydrogel, Advanced Materials, 2007, 19, 1622-1626.
\bibitem{pa9} J.A. Johnson, N.J. Turro, J.T. Koberstein, et al; Some hydrogel having novel molecular structure, Progress in Polymer Science, 2010, 35:332-337.
\bibitem{pa10} S.S. Jang, M. Yashaar, S. Kalani, et al; Mechanical and transport properties of the poly(ethylene oxide)每poly(acrylic acid) double network hydrogel from molecular dynamic simulations, J Phys Chem B 2007;111:1729每37.
\bibitem{pa11} S.S. Jang, M. Yashaar, S. Kalani, et al; Corrections to mechanical and transport properties of the poly(ethylene oxide)每poly(acrylic acid) double network hydrogel from molecular dynamic simulations. J Phys Chem B 2007;111,14440/1.
\bibitem{pa12} K. Luo; The morphology and dynamics of polymerization of polymerization-induced phase separation, European Polymer Journal, 2006, 42: 1499-1505.
\bibitem{pa13} D. Nwabunma, H.W. Chiu, T. Kyu; Theoretical investigation on dynamics of photopolymerization-induced
phase separation and morphology development in nematic liquid crystal/polymer mixtures, J. Chem. Phys (2000). 113, 6429.
\bibitem{pa14} Y. Okumura, K. Ito, The Polyrotaxane Gel: A Topological Gel by Figure-of-Eight Cross-links, Advanced matter. 2001, 13, 485.
\bibitem{pa15} A. Okada, A. Usuki; Twenty years of polymer-clay nanocomposites. Macromol Mater Eng 2006;291:1449每76.
\bibitem{pa16} Z.H. Qiao, Z.R Zhang, T. Tang; An Adaptive time-stepping strategy for the Molecular Beam Epitaxy Models, SIAM J.SCI.COMPUT Vol. 33, No. 3, pp. 1395每1414.
\bibitem{pa17} G. Soderlind; Automatic control and adaptive time-stepping,Numer. Algorighms, 31 (2002),pp. 281每310.
\bibitem{pa18} G. Soderlind, L. Wang; Adaptive Time-Stepping and Computational Stability, Journal of Computational Methods,in Sciences and Engineering,vol. 2, no. 3, 2, pp. 1每3.
\bibitem{pa19} Z. J. Tan, Z. R. Zhang, Y. Q. Huang et al; Moving mesh methods with locally varying time steps, J. Comput. Phys., 200 (2004), pp. 347每367.
\bibitem{pa20} X.C. Xiao, L.Y. Chu, W.M. Chen, et al; Monodispersed thermosponsive hydrogel micropsheres with a volume phase transition driven by hydrogen boning[J]. Ploymer 2005, 46(9):3199-3209.
\bibitem{pa21} Q.H. Zeng, A.B. Yu, G.Q. LU; Multiscale modeling and simulation of polymer nanocomposite, Progress in Polymer Science, 2008, 33: 199-269.
\bibitem{pa22} D. Zhai, H. Zhang; Investigation on the application of the TDGL equation in macromolecular microsphere composite hydrogel, Soft Matter RSC Publishing, 2012. 10. 1039.
\bibitem{pa23} W. Zhang, T.J. Li, P.W. Zhang; Numercial Study for the Nucleation of one-dimensional Stochastic Cahn-Hillard Dynamics, to appear in Commun. Math. Sci., (2012)
\end{thebibliography}
\end{document}